\documentclass{article}
\usepackage{arxiv}
\usepackage[unicode, pdftex]{hyperref}
\usepackage{url}             
\usepackage{booktabs}        
\usepackage{amsfonts}        
\usepackage{nicefrac}        
\usepackage{microtype}       
\usepackage{graphicx}
\usepackage{amsmath}
\usepackage{amssymb}  
\usepackage{natbib}
\usepackage{doi}
\usepackage{xcolor}
\usepackage[all]{xy}
\usepackage{caption}
\usepackage{subcaption} 
\usepackage{orcidlink}
\usepackage{tikz}
\usepackage{amsthm}
\usepackage{natbib}
\usepackage{algorithm}
\usepackage{algpseudocode}  
\usepackage{float}
\newcounter{cdn}
\tikzset{B/.style={blue!50!black},R/.style={red}}
\newcommand{\cd}[3][.5]{\stepcounter{cdn}\begin{tikzpicture}[scale=#1,baseline=-2pt]
\draw circle(1);
\node[anchor=east,font=\small,inner sep=0pt] at (1.22,-1.1){\thecdn};
\foreach\a/\b/\c in{#3}{\draw[\c]({90-360*\a/#2}:1)--({90-360*\b/#2}:1);}
\end{tikzpicture}}

\newtheorem{theorem}{Theorem}[section]
\newtheorem{lemma}[theorem]{Lemma}

\newtheorem{proposition}[theorem]{Proposition}
\theoremstyle{definition}
\newtheorem{definition}[theorem]{Definition}
\theoremstyle{remark}
\newtheorem{remark}[theorem]{Remark}

\title{Enumeration of One-Cycle Signed Chord Diagrams and Their Duals for Optimal Morse Flows on Closed Surfaces}

\renewcommand{\shorttitle}{Enumeration of one-cycle signed chord diagrams}

\author{
  Alexandr Prishlyak \orcidlink{0000-0002-7164-807X}
  \\
  Taras Shevchenko National University of Kyiv, Ukraine
  \\
  \texttt{prishlyak@knu.ua}
}

\hypersetup{
pdftitle={Enumeration of one-cycle signed chord diagrams and their duals for optimal Morse flows on closed surfaces},
pdfauthor={Alexandr Prishlyak},
pdfkeywords={Morse flow, topological classification, vector field, signed chord diagram, dual diagram, enumeration algorithm, surface topology},
colorlinks=true,
linkcolor=blue,
citecolor=blue,
urlcolor=blue
}
\begin{document}
\maketitle

\begin{abstract}
Optimal Morse vector fields on closed surfaces, i.e.\ fields with the minimal possible
number of critical points (a single source and a single sink), are classified up to
topological trajectory equivalence by one-cycle signed chord diagrams
\cite{Kybalko2018}.
 We recall this construction and the associated duality, which
corresponds to reversing the direction of the flow, in a unified encoding in which both
endpoints of a chord carry the same label, and we make it algorithmic. We give an explicit
traversal procedure that detects one-cycle diagrams, a canonical form of an equivalence
class under cyclic shifts, reversal and relabeling of chords, and a construction of the
dual diagram. On this basis we describe an algorithm that enumerates all one-cycle signed
chord diagrams with a given number of chords, both for orientable and for non-orientable
surfaces, groups them into equivalence classes, pairs each class with its dual, and
identifies the self-dual ones. The algorithm is implemented in Python and archived in
Zenodo; its output is validated against the Harer--Zagier numbers in the orientable case
and against orbit-size counts in both cases. We also note that, when a one-cycle diagram is
read as a cellular embedding of a one-vertex graph, its dual diagram describes the
embedding of the dual graph.
\end{abstract}




\section{Introduction}
\label{sec:1}

The topological classification of dynamical systems is a foundational problem in geometric
topology and the qualitative theory of differential equations. Structurally stable systems
are of particular importance because their qualitative phase portraits remain invariant
under small $C^1$-perturbations of the vector field \cite{Peixoto1959, Palis1982}. On
closed surfaces, such structurally stable flows are precisely the Morse--Smale vector fields \cite{Peixoto1973, Oshemkov1998, Kibkalo2022}. Gradient-like
systems and Morse flows arise naturally as gradient fields of Morse functions with respect
to suitable Riemannian metrics \cite{Smale1960, Smale1961}.

According to Smale's construction, an optimal Morse flow, that is, one with a single source
and a single sink, specifies a polar Morse function once all saddle points are placed at the
same critical level. In Fomenko's topological classification of integrable Hamiltonian
systems and of functions on surfaces \cite{Bolsinov_2004}, such a structure is a simple atom
with one input and one output, which corresponds to an edge of degree $2$ in the associated
Reeb (Kronrod--Reeb) graph \cite{Reeb1946, Kronrod1950, Kravchenko_2019, Teresh2024, Tereshchenko2025}.
The global topological structure of such flows is therefore encoded by one-cycle signed
(marked) chord diagrams \cite{Kybalko2018}.

Fine topological invariants and combinatorial structures for flows and functions on surfaces
have been developed in many settings. For closed two-dimensional manifolds, combinatorial
schemes, three-color graphs, chord diagrams and $F$-atoms were used to classify optimal
flows, simple Morse functions and functions with non-degenerate critical points
\cite{prishlyak1997graphs, Prishlyak2000, Prishlyak2002b, Kybalko2018, Prishlyak2016,
Prishlyak2020}. These methods were extended to surfaces with boundary and to non-compact
surfaces, where complete classifications and deformation theories are available
\cite{prishlyak2007classification, Hladysh2016, hladysh2017topology, hladysh2019simple, Hladysh2019, prislyak2017morse, Prishlyak2019,
prishlyak2020three, Borodzik_2016, prishovts26}, and to flows on non-orientable surfaces such
as the projective plane, the M\"obius strip, the Klein bottle and Boy's and Girl's surfaces
\cite{Prishlyak2022, Prishlyak_2023, Prish-Stas2025, PrishButok2025}. In higher dimensions,
topological invariants, handle decompositions, Heegaard diagrams and conjugacy
classifications have been obtained for flows and functions on $3$- and $4$-manifolds
\cite{prish1998, Prishlyak1999, Prishlyak1999a, Prishlyak2001, prishlyak2002morse1, Prishlyak2003, prishlyak2001conjugacy, prishlyak2001, Prishlyak2002beh2,
prishlyak2003regular, prishlyak2005complete, Prishlyak2007, Prishlyak2013, Hatamian2020,
Prishlyak2023, Prishlyak2023a}. Related work on codimension-one flows, bifurcations, closed
$1$-forms
and collective dynamics can be found in \cite{Bilun2002,
Prishlyak2021,  prishovts25}.

Chord diagrams, that is, perfect matchings of points on a circle, are also a classical
object of enumerative combinatorics. The numbers of chord diagrams and of gluings of
polygon sides into surfaces of given genus were studied by Harer and Zagier
\cite{HarerZagier1986}, and an analogue of their formula for unicellular bicolored maps was
obtained in \cite{Adrianov1997}. Counting chord diagrams of a given genus up to rotations was done by
\citet{CoriMarcus1998}, who in particular treated the maximal case, that is, diagrams of
genus $g$ with $2g$ chords; \citet{Krasko2017} extended this count to the full dihedral
group. The orientable one-cycle diagrams of the present paper are exactly these maximal
diagrams. Enumeration of chord diagrams of other types, motivated by Vassiliev knot
invariants, and of diagrams up to rotation and reflection is given in
\cite{Stoimenow2000, Khruzin2000}. Chord diagrams were
further used as an invariant of functions on surfaces by Kadubovsky, who classified
functions with one essentially critical point on oriented surfaces with boundary by means of
color chord diagrams \cite{Kadubovsky2003}, studied the realization and topological
equivalence of functions with a single saddle point on oriented surfaces and counted
non-equivalent functions for surfaces of prime-related genus \cite{Kadubovskyi2006}, and
computed the numbers of non-equivalent minimal functions on closed oriented and
non-orientable surfaces \cite{Kadubovskyi2015}. In a series of papers by Kadubovskyi and
his co-authors, two-color chord diagrams with a prescribed number of faces are enumerated
under the action of the cyclic and the dihedral groups, with explicit formulas and tables of
values \cite{KadubovskyiSypchuk2020, KadubovskyiStepkinKyrychenko2021}. Our one-cycle
diagrams belong to the same circle of enumeration problems, but we treat the orientable and
the non-orientable case, all genera and the pairing with duals within a single algorithm.

The chord diagram construction for optimal flows and the duality that corresponds to time
reversal are due to \citet{Kybalko2018}. The aim of the present paper is to make this
theory algorithmic. In contrast to the classical notation, in which the two endpoints of a
chord carry different indices, each chord here receives one label $i\in\{1,\dots,k\}$ that
occurs twice along the circle, which gives a canonical chord code of length $2k$. We recall
the connection between one-cycle chord diagrams, Reeb graphs and embeddings of one-vertex
graphs, and we give an explicit cycle-tracing procedure, a canonical form of an equivalence
class, a construction of the dual diagram, and an enumeration algorithm for all one-cycle
signed chord diagrams with a given number of chords. The algorithms are implemented in
Python and archived together with visualization tools on Zenodo \cite{zenodo22298943,
zenodo23084341}. 

We describe an algorithm that, for given $g$, enumerates all one-cycle chord diagrams and identifies the pairs of dual diagrams, corresponding to the flows $X$ and $-X$, together with the self-dual diagrams, whose flows are equivalent to their reverses.

The paper is organized as follows. Section~\ref{sec:prelim} collects the necessary facts on
Morse--Smale flows, optimal Morse flows and topological equivalence. Section~\ref{sec:chord}
introduces signed chord diagrams, the code with identical endpoint labels and the canonical
form of an equivalence class. Section~\ref{sec:construction} recalls the construction of the
chord diagram of an optimal flow and shows that it is a one-cycle diagram.
Section~\ref{sec:dual} treats the chord cycle and the dual diagram, and
Section~\ref{sec:software} describes the enumeration algorithm, its implementation and the
numerical results.


\section{Preliminaries}\label{sec:prelim}

\subsection{Morse--Smale and optimal Morse flows}\label{subsec:ms}

We consider $C^r$-smooth ($r\ge1$) dynamical systems with continuous time (flows) on a
smooth compact surface $M$. Let $X$ be a vector field generating the flow $X_t$; the
trajectory through $x\in M$ is denoted by $\varphi_t(x)=X_t(x)$.

\begin{definition}
Fixed (singular) points and closed trajectories are called the \textit{critical elements}
of $X$.
\end{definition}

\begin{definition}
A singular point $p$ of $X$ is \textit{hyperbolic} (non-degenerate) if the Jacobian matrix
\begin{equation}\label{eq:jacobian}
\left( \frac{\partial X_i}{\partial x_j}(p) \right)_{i,j=1}^{2},
\end{equation}
written in local coordinates, has no eigenvalue with zero real part. A closed trajectory
$\gamma$ is \textit{hyperbolic} if the differential at $p\in\gamma$ of the Poincar\'e return
map has no eigenvalue of absolute value $1$.
\end{definition}

\begin{definition}
A point $y\in M$ is an $\omega$-limit ($\alpha$-limit) point of $x\in M$ if there is a
sequence $t_i\to+\infty$ ($t_i\to-\infty$) such that $\varphi_{t_i}(x)\to y$. The set of all
such points is the $\omega$-limit ($\alpha$-limit) set $\omega(x)$ ($\alpha(x)$).
\end{definition}

\begin{definition}
The \textit{stable} and \textit{unstable manifolds} of a critical element $c$ are
\[
W^s(c)=\{x\in M \mid \omega(x)=c\},\qquad W^u(c)=\{x\in M \mid \alpha(x)=c\}.
\]
For a singular point $p$ this means $\lim_{t\to+\infty}\varphi_t(x)=p$ and
$\lim_{t\to-\infty}\varphi_t(x)=p$, respectively; in particular $p\in W^s(p)\cap W^u(p)$.
\end{definition}

\begin{definition}
A point $p\in M$ is \textit{wandering} for $X$ if there are a neighborhood $V$ of $p$ and
$T>0$ such that $X_t(V)\cap V=\emptyset$ for all $|t|>T$; otherwise $p$ is
\textit{non-wandering}. The set of non-wandering points is denoted by $\Sigma(X)$.
\end{definition}

\begin{definition}
A vector field $X$ on a compact surface $M$ is a \textit{Morse--Smale vector field} if
(1) $X$ has finitely many critical elements, all hyperbolic; (2) $\Sigma(X)$ is the union of
the critical elements; (3) the stable and unstable manifolds of critical elements intersect
transversally.
\end{definition}

Hyperbolic singular points of a flow on a surface are of three types: \textit{sinks},
\textit{saddles} and \textit{sources}. With the convention $X=-\operatorname{grad} f$ for a
Morse function $f$, they correspond to critical points of index $0$, $1$ and $2$,
respectively. Transversality excludes trajectories going from one saddle to another
(saddle connections), since such a trajectory would lie in both $W^u$ and $W^s$ with the same
tangent direction.

\begin{definition}
A Morse--Smale flow is a \textit{Morse flow} if it has no closed trajectories. A Morse flow
on a closed surface $M$ is \textit{optimal} if it has the minimal number of singular points
among all Morse flows on $M$.
\end{definition}

\begin{theorem}\label{thm:optimal}
Let $M$ be a closed connected surface. Then:
\begin{enumerate}
\item[(a)] a Morse flow on $M$ is optimal if and only if it has exactly one sink and exactly
  one source;
\item[(b)] an optimal Morse flow has $k=2-\chi(M)$ saddles; thus $k=2g$ if $M$ is orientable
  of genus $g$, and $k=h$ if $M$ is non-orientable of genus $h$.
\end{enumerate}
Optimal Morse flows exist on every closed connected surface.
\end{theorem}
\begin{proof}
Let $a$, $b$, $k$ be the numbers of sources, sinks and saddles of a Morse flow. The sum of
the indices of the singular points is $\chi(M)$ (Poincar\'e--Hopf), and a source or sink has
index $+1$, a saddle $-1$; hence $a+b-k=\chi(M)$ and the total number of singular points is
$a+b+k=2(a+b)-\chi(M)$. Since $\Sigma(X)$ is the finite set of singular points, $\omega(x)$ is
a single singular point $p$ for every $x$. Take $x$ outside the one-dimensional stable
manifolds of the saddles and outside the singular points; then $p$ is not a saddle, and not a
source (as $W^s(p)=\{p\}$ for a source), hence a sink. So $b\ge1$, and by time reversal $a\ge1$.
Therefore $a+b\ge2$ with equality exactly when $a=b=1$, which gives (a) and $k=2-\chi(M)$.
Existence is shown in \cite{Kybalko2018} (equivalently, via a Morse function with one
maximum and one minimum).
\end{proof}

\subsection{Topological equivalence}\label{subsec:equiv}

\begin{definition}\label{def:equiv}
Two vector fields $X,Y$ on $M$ are \textit{topologically (trajectory) equivalent},
$X\sim Y$, if there is a homeomorphism $h\colon M\to M$ that maps trajectories of $X$ onto
trajectories of $Y$ and preserves the direction of motion along them. The homeomorphism $h$
is not required to preserve an orientation of $M$.
\end{definition}

The time parametrization is ignored. The \textit{inverse} of $X$ is the field $-X$, with the
same trajectories traversed in the opposite direction. In general $X\not\sim -X$, so the
relation ``$X\sim -Y$'' is a different, weaker condition; both occur in Section~\ref{sec:dual}.

\begin{theorem}
Let $h$ be a topological equivalence of $X$ and $Y$. Then $p$ is a singular point of $X$ iff
$h(p)$ is a singular point of $Y$; a trajectory of $X$ is closed iff its image is closed; and
$h$ maps the $\alpha$-limit ($\omega$-limit) set of every trajectory onto the $\alpha$-limit
($\omega$-limit) set of the image trajectory.
\end{theorem}

\begin{definition}
A vector field $X\in\mathcal X^r(M)$, $r\ge1$, is \textit{structurally stable} if it has a
neighborhood $U(X)$ in the $C^r$-topology such that $Y\sim X$ for all $Y\in U(X)$.
\end{definition}

On closed surfaces the structurally stable fields are exactly the Morse--Smale fields
\cite{Peixoto1973, Oshemkov1998, Kibkalo2022}; in particular, optimal Morse flows are
structurally stable, and their classification up to $\sim$ is the object of this paper.


\section{Signed Chord Diagrams and Canonical Encoding}\label{sec:chord}

We formalize chord diagrams using a unified endpoint labeling system suitable for
algorithmization.

\begin{definition}
A \textit{chord diagram} $D$ with $k$ chords consists of a boundary circle $S^1$,
$2k$ distinct points on $S^1$, and $k$ chords connecting these points in pairs
(chords may cross when drawn in the disk).
\end{definition}

In our encoding scheme:
\begin{itemize}
    \item Each chord is assigned a unique identifier $i \in \{1, 2, \dots, k\}$.
    \item Both endpoints of chord $i$ are assigned the same integer label $i$.
    \item Traversing $S^1$ counterclockwise from a fixed origin yields a sequence
          $\mathbf{c} = (c_0, c_1, \dots, c_{2k-1})$ of length $2k$ in which each integer
          from $1$ to $k$ appears exactly twice; positions are indexed by $\mathbb{Z}_{2k}$.
\end{itemize}

\begin{definition}
A code $\mathbf{c}$ is in \textit{Restricted Growth String (RGS) canonical form} if
$c_0 = 1$, and for every $m$ with $1\le m<k$ the first occurrence of $m$ precedes the
first occurrence of $m+1$. For an arbitrary sequence $\mathbf{w}$ in which every label
occurs twice, its \textit{RGS normalization} relabels the chords in order of first
occurrence.
\end{definition}

The RGS codes of length $2k$ are in bijection with the perfect matchings of $2k$ points;
there are $(2k-1)!!$ of them.

\begin{definition}
A \textit{signed chord diagram} (or marked chord diagram) assigns a sign
$\sigma_i \in \{+1, -1\}$ to each chord $i \in \{1, \dots, k\}$.
\end{definition}

\begin{definition}
Two signed chord diagrams $D_1, D_2$ are \textit{isomorphic} if one can be transformed
into the other via a cyclic shift of the code and a relabeling of the chords (signs are
carried along with their chords). They are \textit{equivalent} if they can additionally
be related by a reversal (mirror reflection) of the boundary circle. We write
$D_1\simeq D_2$ for equivalence.
\end{definition}

\subsection{Canonical form of an equivalence class}\label{subsec:canon}

Let $\mathcal G$ be the dihedral group of order $4k$ generated by the cyclic shift and
the reversal of the code. For $g\in\mathcal G$ the diagram $g\cdot D$ is obtained by
applying $g$ to the code, taking the RGS normalization, and carrying the signs along with
their chords (shifts and reversals do not change signs; relabeling permutes them). The
\textit{canonical form} of $D$ is the lexicographic minimum
\begin{equation}\label{eq:canon}
  \operatorname{can}(D)=\min_{g\in\mathcal G}\ g\cdot D .
\end{equation}
Two signed diagrams are equivalent iff their canonical forms coincide. It is computed in
$O(k^2)$ time ($4k$ images, each normalized in $O(k)$). If $\mathcal G$ is replaced by
its cyclic subgroup of order $2k$ (shifts only), the same procedure yields isomorphism
classes. The orbit of $D$ consists of $4k/|\operatorname{Stab}(D)|$ distinct normalized
codes; the stabilizer is the group of symmetries of $D$.


\section{Construction of a Chord Diagram from an Optimal Flow}\label{sec:construction}

The construction of a chord diagram from an optimal Morse flow follows
\citet{Kybalko2018}; we recall it here in the notation of Section~4.

Let $M$ be a closed connected surface (orientable or not) and $X$ an optimal Morse vector
field on $M$, i.e.\ a Morse flow with a single source $x_0$, a single sink $y_0$, and
$k$ saddle points $s_1,\dots,s_k$. By Theorem~\ref{thm:optimal}, the number of saddles is $k=2-\chi(M)$, that is, $k=2g$ if $M$ is orientable of genus $g$ and $k$ equals the non-orientable genus otherwise.

\textbf{Neighborhood of the source.}
Let $U(x_0)$ be a small closed disk around $x_0$ whose boundary $\partial U(x_0)\cong S^1$
is transversal to the trajectories of $X$ leaving $x_0$; each such trajectory meets
$\partial U(x_0)$ exactly once.

\textbf{Chords.}
For each saddle $s_i$, the set $W^s(s_i)\setminus\{s_i\}$ consists of two trajectories
(separatrices) issuing from $x_0$ and entering $s_i$. Hence $\overline{W^s(s_i)}$ meets
$\partial U(x_0)$ in exactly two points. Stable manifolds of distinct saddles are
disjoint, so these $k$ pairs of points are pairwise disjoint on $\partial U(x_0)$, which
gives $2k$ distinct points. Joining the two points of each pair by a chord inside the
circle (the chords may cross each other in the picture) and labeling it by $i$ yields a
chord diagram with $k$ chords.

\textbf{Signs.}
Let $V_i$ be a regular neighborhood of $U(x_0)\cup\overline{W^s(s_i)}$, i.e.\ the disk
$U(x_0)$ with one band attached along the two points of the chord $i$. If $V_i$ is
orientable (an annulus), put $\sigma_i=+1$; if $V_i$ is a M\"obius band, put
$\sigma_i=-1$. (In terms of the marking, $+1$ corresponds to an unmarked and $-1$ to a
marked chord.)

\begin{proposition}\label{prop:one_cycle}
The signed chord diagram $D_X$ constructed above is a one-cycle diagram. Its chord cycle
is the boundary of the basin of attraction of the sink $y_0$.
\end{proposition}
\begin{proof}
The set $\Gamma=\{x_0\}\cup\bigcup_{i}\overline{W^s(s_i)}$ is a bouquet of $k$ circles
embedded in $M$ (the loop number $i$ is formed by the two separatrices entering $s_i$).
Every point of $M\setminus\Gamma$ lies on a trajectory with $\omega$-limit $y_0$, so
$M\setminus\Gamma$ is the basin of attraction of the sink, which is an open disk.
Hence $\Gamma$ is a cellular embedding with a single face, and the boundary of this face,
read along $\partial U(x_0)$ and the chords, is exactly the traversal
\eqref{eq:step} with the signs defined above.
\end{proof}

\textbf{Dependence on choices.}
The diagram depends on the choice of the disk $U(x_0)$, of the initial point on
$\partial U(x_0)$, of the orientation of $\partial U(x_0)$ (we do not fix an orientation
of $M$), and of the numbering of the saddles. Changing them results, respectively, in an
isotopy of the diagram, a cyclic shift of the code, a reversal, and a relabeling of the
chords; the signs are unchanged by the first three operations and are carried along with
their chords by the last one. Thus every optimal Morse flow on a closed surface induces a
well-defined equivalence class of signed one-cycle chord diagrams (Section~\ref{sec:chord}).
Conversely, every equivalence class of one-cycle diagrams is realized by an optimal flow \citep[Theorem~3]{Kybalko2018}.

\section{Dual Chord Diagrams and Trajectory Cycles}\label{sec:dual}

We formalize the dualization procedure (formerly termed inverse diagrams) under the
identical-label code format.

\subsection{Tracing the chord cycle}\label{subsec:trace}

Informally: start on a boundary arc between two chord endpoints and move
counterclockwise; on reaching an endpoint of chord $i$, jump along the chord to its
opposite endpoint; then continue along $S^1$ in the same direction if $\sigma_i=+1$, and
in the reversed direction if $\sigma_i=-1$; repeat until the starting arc is reached again
in the same direction.

Formally, for a position $p$ let $\pi(p)$ be the position of the second endpoint of the
chord $c_p$. A \textit{state} is a pair $(p,\varepsilon)$ with $p\in\mathbb Z_{2k}$ and
$\varepsilon\in\{+1,-1\}$: the walker stands at the point $p$ and is about to move along
the boundary arc towards the neighbouring point $p+\varepsilon$. One step is
\begin{equation}\label{eq:step}
  (p,\varepsilon)\ \longmapsto\ \bigl(\pi(q),\ \varepsilon\,\sigma_{c_q}\bigr),
  \qquad q=p+\varepsilon ,
\end{equation}
where $q$ is called the \textit{arrival point}. The tracing starts in the state $(0,+1)$,
i.e.\ on the arc right after the first element of the code, and stops when this state
recurs. The arrival points $q_1,\dots,q_m$ define the \textit{visit sequence}
$\mathbf w=(c_{q_1},\dots,c_{q_m})$.

\begin{algorithm}[h]
\caption{Tracing the chord cycle}\label{alg:trace}
\begin{algorithmic}[1]
\Require code $\mathbf c$, signs $\sigma$
\State $(p,\varepsilon)\gets(0,+1)$; \ $Q\gets[\,]$ \Comment{list of arrival points}
\Repeat
  \State $q\gets (p+\varepsilon)\bmod 2k$; \ append $q$ to $Q$
  \State $\varepsilon\gets \varepsilon\cdot\sigma_{c_q}$; \ $p\gets\pi(q)$
\Until{$(p,\varepsilon)=(0,+1)$}
\State \Return $Q$
\end{algorithmic}
\end{algorithm}

The map \eqref{eq:step} is a permutation of the $4k$ states, so the loop terminates in at
most $4k$ steps.

\begin{definition}
A \textit{chord cycle} is the closed traversal path generated by \eqref{eq:step}. A
diagram is a \textit{one-cycle diagram} if $|Q|=2k$ and every label $1,\dots,k$ occurs
exactly twice in $\mathbf w$; equivalently, a single traversal covers all chords twice and
all $2k$ boundary arcs.
\end{definition}

For $\sigma\equiv+1$ the cycle passes every chord once in each direction. For a chord
with $\sigma_i=-1$ both passes may arrive at the same endpoint, so the arrival points
need not be pairwise distinct. For a one-vertex one-face map with $k$ edges the Euler
characteristic is $2-k$; hence one-cycle diagrams with all signs $+1$ exist only for
$k=2g$ (orientable surface of genus $g$), and one-cycle diagrams with at least one sign
$-1$ describe non-orientable surfaces of genus $k$.

\subsection{Dual chord diagram}\label{subsec:dual}

Orient each chord $i$ from its endpoint $a_i$ to its endpoint $b_i$, where $a_i<b_i$ are
the positions of its endpoints in $\mathbf c$. In the visit sequence $\mathbf w$ the label
$i$ occurs twice, and each occurrence is an arrival either at $a_i$ or at $b_i$; the
jump after an arrival at $a_i$ runs along the chord in the positive direction, after an
arrival at $b_i$ in the negative one.

To construct the \textit{dual chord diagram} $D^*$ of a one-cycle diagram $D$:
\begin{enumerate}
\item Trace the chord cycle of $D$ and form the visit sequence $\mathbf w$ of length $2k$.
\item The dual code is the RGS normalization $\mathbf c^*$ of $\mathbf w$; let $\mu$ be
      the relabeling that sends the labels of $D$ to the labels of $\mathbf c^*$.
\item The sign $\sigma^*_{\mu(i)}$ is $+1$ if the cycle runs along chord $i$ in opposite
      directions (the two arrivals are at different endpoints), and $-1$ if it runs along
      $i$ in the same direction (both arrivals at the same endpoint).
\item Replace the result by its canonical form \eqref{eq:canon}.
\end{enumerate}

Since $\mathbf w$ depends on the base point and the direction of the tracing, the
dual is a well-defined \textit{equivalence class}.

\begin{lemma}\label{lem:dual}
If $D$ is a one-cycle diagram, then $D^*$ is a one-cycle diagram, $D^{**}\simeq D$,
$D$ and $D^*$ have the same (non-)orientable genus, and the class of $D^*$ depends only
on the class of $D$. Consequently, if $N$ is the number of equivalence classes,
$K$ the number of unordered pairs $\{[D],[D^*]\}$ and $s$ the number of self-dual
classes, then $s=2K-N$.
\end{lemma}
\begin{proof}
See \cite{Kybalko2018}.
\end{proof}

\begin{theorem}\label{thm:dual_flow}
Two optimal Morse flows $X,Y$ are topologically equivalent if and only if $D_X\simeq D_Y$, and $X$ is topologically equivalent to the inverse flow $-Y$ if and only if $D_X\simeq D_Y^*$.
\end{theorem}
\begin{proof}
See \cite{Kybalko2018}.
\end{proof}

\begin{remark}\label{rem:dual_graph}
If a one-cycle signed chord diagram $D$ is interpreted as a cellular embedding of a
single-vertex graph (a bouquet of $k$ circles) on a closed surface $M$ with a single
2-cell, then the dual chord diagram $D^*$ specifies the cellular embedding of the
topological dual graph $G^*$ on $M$; $G^*$ is again a bouquet of $k$ circles with a
single 2-cell.
\end{remark}


\section{Software Implementation and Enumeration Data}\label{sec:software}

The enumeration algorithm has been fully implemented in Python and deposited in Zenodo
\cite{zenodo22298943}. For the genus parameter $g$ the program takes $k=2g$ and
$\sigma\equiv(+1)^k$ if $g\ge0$ (orientable case), and $k=-g$ with all $2^k-1$ sign
vectors different from $(+1)^k$ if $g<0$ (non-orientable case).

\begin{algorithm}[h]
\caption{Enumeration of one-cycle classes}\label{alg:enum}
\begin{algorithmic}[1]
\Require number of chords $k$, orientability flag
\State $\mathcal C\gets\emptyset$
\For{each RGS code $\mathbf c$ of length $2k$} \Comment{$(2k-1)!!$ codes}
  \For{each admissible sign vector $\sigma$}
    \If{$(\mathbf c,\sigma)$ is a one-cycle diagram}
      \State add $\operatorname{can}(\mathbf c,\sigma)$ to $\mathcal C$
    \EndIf
  \EndFor
\EndFor
\State $\mathcal P\gets$ the set of unordered pairs $\{C,\,C^*\}$ with $C\in\mathcal C$
\State \Return $(\mathcal C,\mathcal P)$ \Comment{a pair with $C=C^*$ is self-dual}
\end{algorithmic}
\end{algorithm}

RGS codes are generated by backtracking: at each position one either closes one of the
currently open chords or opens the next unused label, so exactly the $(2k-1)!!$ perfect
matchings of the $2k$ points appear, each once. Each code is tested by Algorithm
\ref{alg:trace}; the total cost of the tests is $O\bigl((2k-1)!!\cdot 2^{k}\cdot k\bigr)$
(the factor $2^k$ is absent in the orientable case), plus $O(k^2)$ for each code that
passes. The set $\mathcal C$ is a hash set, so the memory is proportional to $N$. For
every class the dual is computed as in Subsection \ref{subsec:dual} and located in
$\mathcal C$ (Lemma \ref{lem:dual} guarantees that it is there). The program outputs
$N=|\mathcal C|$, the number $K=|\mathcal P|$ of pairs and the number $s=2K-N$ of
self-dual classes. The list of classes can be saved to a JSON file containing the number
of chords, the input genus, the orientability flag, and a list of diagrams, each given by
its canonical code and its sign vector (one representative per class); the pairs are
recovered by applying the duality of Subsection \ref{subsec:dual}.

\paragraph{Validation.}
The program was checked against known enumerative results. In the orientable case the
number of one-cycle \emph{rooted} codes equals the Harer--Zagier number
$\varepsilon_g(2g)=\frac{(4g)!}{4^g(2g+1)!}$ \cite{HarerZagier1986}: $1,\,21,\,1485$ for
$g=1,2,3$. The numbers $N=1,\,4,\,82$ of equivalence classes under the dihedral group for
$g=1,2,3$ coincide with the numbers of maximal chord diagrams up to all isomorphisms
computed by \citet{Krasko2017}. For arbitrary sign vectors, the numbers of rooted one-cycle
codes (orientable and non-orientable together) $1,\,5,\,41,\,509,\,8229$ for $k=1,\dots,5$
coincide with the numbers of rooted one-vertex one-face maps with $k$ edges on all closed
surfaces given there. Moreover, in both cases the sum of the orbit sizes
$4k/|\operatorname{Stab}|$ over all classes equals the number of rooted one-cycle codes
(e.g.\ $21$ for $k=4$ orientable and $41$ for $k=3$ non-orientable), and for every class its
dual is found in $\mathcal C$.

\begin{table}[ht!]
\centering
\caption{Enumeration metrics for one-cycle signed chord diagrams across surface genus $g$. Here $N$ is the total count of non-equivalent diagrams, $K$ is the number of distinct dual pairs, and $2K-N$ is the number of self-dual diagrams.}
\label{tab:results}
\begin{tabular}{|c|c|c|c|}
\hline
\textbf{Genus $g$} & \textbf{Total $N$} & \textbf{Dual Pairs $K$} & \textbf{Self-Dual $2K - N$} \\ \hline
-8 & 3534490 & 1775258 & 16026 \\
-7 & 144904 & 73220 & 1536 \\
-6 & 7190 & 3847 & 504 \\
-5 & 473 & 266 & 59 \\
-4 & 47 & 35 & 23 \\
-3 & 8 & 6 & 4 \\
-2 & 2 & 2 & 2 \\
-1 & 1 & 1 & 1 \\
0  & 1 & 1 & 1 \\
1  & 1 & 1 & 1 \\
2  & 4 & 4 & 4 \\
3  & 82 & 56 & 30 \\
4  & 7258 & 3827 & 396 \\ \hline
\end{tabular}
\end{table}


\section{Conclusion}

In this paper, we have established a unified combinatorial framework for the topological classification of optimal Morse flows on closed surfaces based on single-cycle chord diagrams. By leveraging Smale's theorem, we demonstrated that setting equal values at all saddle critical points naturally translates polar Morse flows into Morse functions corresponding to Fomenko atoms with a single input and output, which map to degree-$2$ edges on Kronrod--Reeb graphs. The upgraded indexing convention--assigning matching labels $i \in \{1, \dots, k\}$ to chord endpoints--provides a compact and canonical representation suitable for algorithmic enumeration and dual graph embeddings.

The upgraded indexing convention --- assigning matching labels $i\in\{1,\dots,k\}$ to chord endpoints --- provides a compact and canonical representation suitable for algorithmic enumeration. Moreover, the dual chord diagram describes the same flow with reversed time: replacing the vector field $X$ by $-X$ interchanges the source and the sink, and the diagram is replaced by its dual. Our algorithm therefore not only lists the diagrams but also finds the pairs of dual diagrams, i.e.\ the pairs of mutually reversed flows, and the self-dual diagrams, i.e.\ the flows equivalent to their reverses. For example, among the $82$ diagrams of genus $g=3$ there are $30$ self-dual ones and $26$ dual pairs, so there are $56$ classes of flows up to time reversal.

The methodology developed here opens several natural avenues for further research and generalizations:
\begin{enumerate}
    \item \textbf{Multi-cycle chord diagrams:} Extending the canonical chord encoding from single-source (polar) flows to arbitrary Morse-Smale flows with multiple sources and sinks on higher-genus closed surfaces.
    \item \textbf{Manifolds with boundary:} Adapting the chord-diagram formalism to capture boundary critical points, sliding vector fields, and boundary-parallel trajectories on surfaces with non-empty boundary.
    \item \textbf{Higher-dimensional generalizations:} Applying multi-dimensional handle decomposition codes and generalized Reeb graphs to classify Morse-Smale flows on 3-manifolds and 4-manifolds.
    \item \textbf{Software ecosystem and visualization:} Integrating our single-cycle chord diagram algorithms with interactive software packages on Zenodo \cite{zenodo22298943, zenodo23084341} to facilitate the automated generation, dual transformation, and visual rendering of complex dynamical phase portraits.
\end{enumerate}

\newpage
\section{Appendix. Catalog of One-Cycle Chord Diagrams}

This section contains a comprehensive graphical atlas of representative one-cycle signed chord diagrams generated for $g\in\{-5,-4,\ldots,3\}$. New figures can be inserted into the template structures below.

\subsection{From the computational code to the drawing code}

For computations a one-cycle chord diagram with $k$ chords is stored as a word
$w=(w_1,\dots,w_{2k})$ in which every label $m\in\{1,\dots,k\}$ occurs exactly twice
(the two ends of chord $m$), together with a sign $\varepsilon_m\in\{+1,-1\}$ for each chord.
The word is normalised by first occurrence, so the labels first appear in the order $1,2,\dots,k$.
Orientable diagrams ($g\ge 0$, $k=2g$) have all $\varepsilon_m=+1$; non-orientable ones ($g<0$, $k=-g$) have at least one $\varepsilon_m=-1$.

The drawing code used by the macro \verb|\cd| is a list of triples \verb|a/b/c|, where
$a<b$ are the cyclic numbers of the ends of a chord on the circle with $n=2k$ points,
and \verb|c| is its colour. The conversion is:
\begin{enumerate}
\item Number the ends cyclically $1,\dots,n$ in the order of the word; position $i$ is drawn at the polar angle $90^\circ-360^\circ i/n$, i.e.\ clockwise starting just after the top point.
\item For each label $m$ let $a_m<b_m$ be the two positions with $w_{a_m}=w_{b_m}=m$.
\item Put $c_m=\texttt{B}$ (blue) if $\varepsilon_m=+1$ and $c_m=\texttt{R}$ (red) if $\varepsilon_m=-1$.
\item Output the triples $a_m/b_m/c_m$ for $m=1,\dots,k$ separated by commas. Because of the normalisation $a_1<a_2<\dots<a_k$, no extra sorting is needed.
\end{enumerate}
The diagram is then inserted as \verb|\cd[s]{n}{list}|, with $n=2k$ and scale $s$ (we use $0.62$--$0.65$).
For example, the word $1\,2\,1\,2\,3\,4\,3\,4$ (all signs $+1$) gives
\verb|\cd[0.62]{8}{1/3/B,2/4/B,5/7/B,6/8/B}|.
The empty diagram of the sphere is drawn as \verb|\cd[0.62]{2}{}|.

\subsection{Order, numbering, duality}

Within each genus the diagrams are the representatives of the equivalence classes under rotation, reflection and relabelling, one per class: the representative is the lexicographically smallest normalised word in its class, and the classes are listed in increasing order of this word. The number under a diagram is the counter \verb|cdn|, reset at the start of each genus and continued across figures if the list of one genus does not fit into one figure.

The dual diagram is obtained by traversing the unique face of the diagram: starting at an arc of the circle, we move along the circle to the next chord end, cross the chord to its other end, and continue along the circle in the same direction if the chord is untwisted and in the opposite direction if it is twisted, until we return to the starting arc. The sequence of chords crossed along the way is the word of the dual diagram. The dual chord $m$ is untwisted if the two crossings of chord $m$ go in opposite directions (from its smaller end to its larger one and back), and twisted if both go the same way. The dual is then reduced to its canonical form and found in the same list; a diagram is self-dual if it coincides with its own dual, otherwise it forms a dual pair. After each list we state the numbers of the self-dual diagrams and of the dual pairs.

\subsection{Orientable Surfaces}

\begin{figure}[ht!]
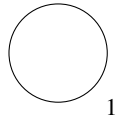

\begin{center}
\setcounter{cdn}{0}
\cd[0.65]{2}{}
\end{center}
\caption{One-cycle chord diagrams for the sphere ($g=0$, $k=0$).}
\label{fig:atlas_g0_1}
\end{figure}


\begin{figure}[ht!]
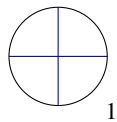

\begin{center}
\setcounter{cdn}{0}
\cd[0.65]{4}{1/3/B,2/4/B}
\end{center}
\caption{One-cycle chord diagrams for the torus ($g=1$, $k=2$).}
\label{fig:atlas_g1_1}
\end{figure}


\begin{figure}[ht!]
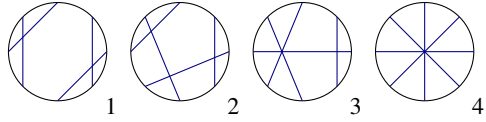

\begin{center}
\setcounter{cdn}{0}
\cd[0.65]{8}{1/3/B,2/4/B,5/7/B,6/8/B}\enskip
\cd[0.65]{8}{1/3/B,2/5/B,4/7/B,6/8/B}\enskip
\cd[0.65]{8}{1/3/B,2/6/B,4/7/B,5/8/B}\enskip
\cd[0.65]{8}{1/5/B,2/6/B,3/7/B,4/8/B}
\end{center}
\caption{One-cycle chord diagrams for genus $g=2$ ($k=4$). 4 diagrams in total: 4 self-dual, 0 dual pairs.
\emph{Self-dual:} 1, 2, 3, 4.}
\label{fig:atlas_g2_1}
\end{figure}


\begin{figure}[ht!]
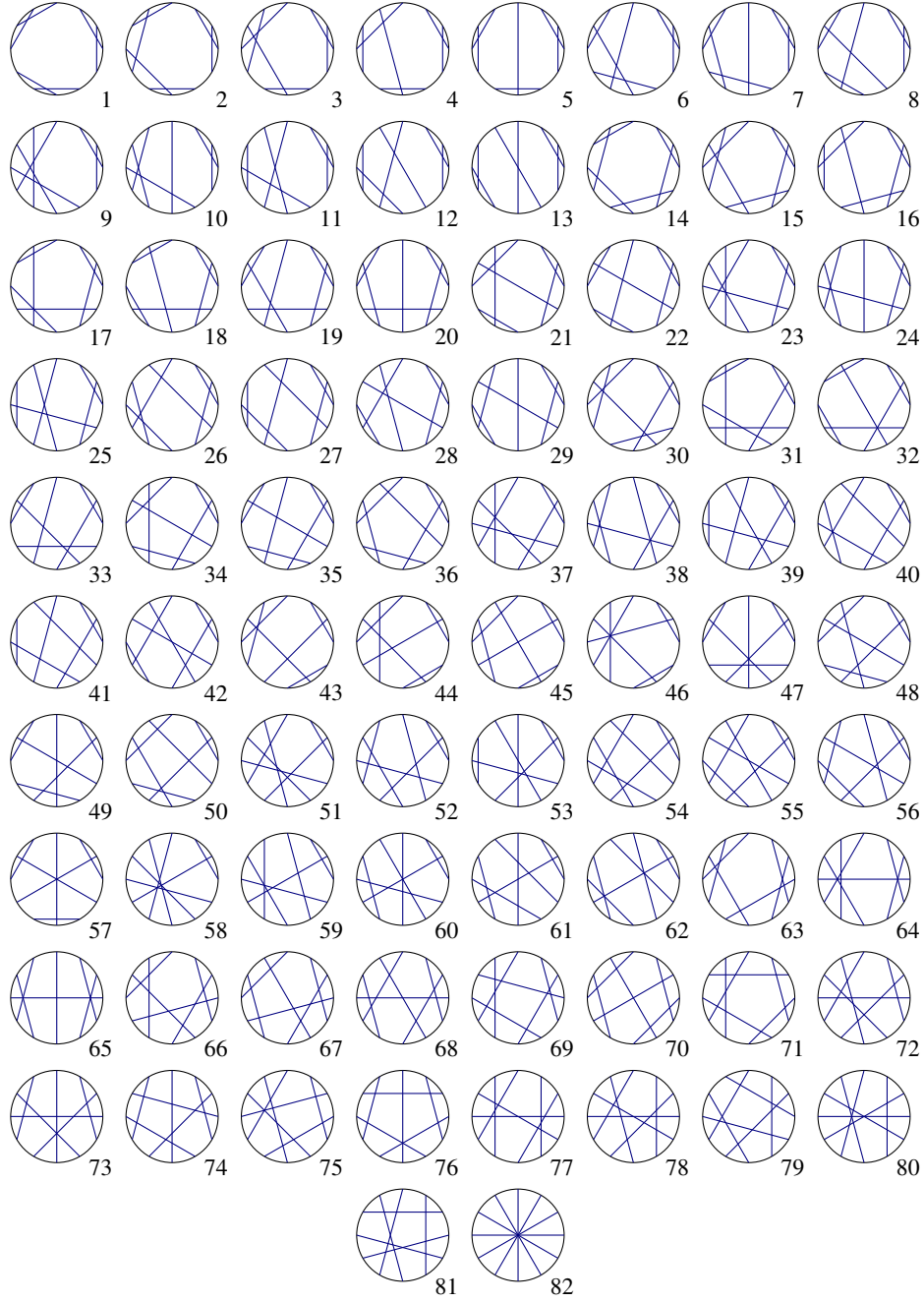

\begin{center}
\setcounter{cdn}{0}
\cd[0.62]{12}{1/3/B,2/4/B,5/7/B,6/8/B,9/11/B,10/12/B}\enskip
\cd[0.62]{12}{1/3/B,2/4/B,5/7/B,6/9/B,8/11/B,10/12/B}\enskip
\cd[0.62]{12}{1/3/B,2/4/B,5/7/B,6/10/B,8/11/B,9/12/B}\enskip
\cd[0.62]{12}{1/3/B,2/4/B,5/7/B,6/11/B,8/10/B,9/12/B}\enskip
\cd[0.62]{12}{1/3/B,2/4/B,5/7/B,6/12/B,8/10/B,9/11/B}\enskip
\cd[0.62]{12}{1/3/B,2/4/B,5/8/B,6/10/B,7/12/B,9/11/B}\enskip
\cd[0.62]{12}{1/3/B,2/4/B,5/8/B,6/12/B,7/10/B,9/11/B}\enskip
\cd[0.62]{12}{1/3/B,2/4/B,5/10/B,6/8/B,7/12/B,9/11/B}\\[4pt]
\cd[0.62]{12}{1/3/B,2/4/B,5/9/B,6/10/B,7/11/B,8/12/B}\enskip
\cd[0.62]{12}{1/3/B,2/4/B,5/9/B,6/12/B,7/10/B,8/11/B}\enskip
\cd[0.62]{12}{1/3/B,2/4/B,5/9/B,6/11/B,7/12/B,8/10/B}\enskip
\cd[0.62]{12}{1/3/B,2/4/B,5/11/B,6/9/B,7/12/B,8/10/B}\enskip
\cd[0.62]{12}{1/3/B,2/4/B,5/11/B,6/12/B,7/9/B,8/10/B}\enskip
\cd[0.62]{12}{1/3/B,2/5/B,4/7/B,6/9/B,8/11/B,10/12/B}\enskip
\cd[0.62]{12}{1/3/B,2/5/B,4/7/B,6/10/B,8/11/B,9/12/B}\enskip
\cd[0.62]{12}{1/3/B,2/5/B,4/7/B,6/11/B,8/10/B,9/12/B}\\[4pt]
\cd[0.62]{12}{1/3/B,2/5/B,4/8/B,6/9/B,7/11/B,10/12/B}\enskip
\cd[0.62]{12}{1/3/B,2/5/B,4/8/B,6/11/B,7/9/B,10/12/B}\enskip
\cd[0.62]{12}{1/3/B,2/5/B,4/8/B,6/10/B,7/12/B,9/11/B}\enskip
\cd[0.62]{12}{1/3/B,2/5/B,4/8/B,6/12/B,7/10/B,9/11/B}\enskip
\cd[0.62]{12}{1/3/B,2/5/B,4/10/B,6/8/B,7/11/B,9/12/B}\enskip
\cd[0.62]{12}{1/3/B,2/5/B,4/10/B,6/8/B,7/12/B,9/11/B}\enskip
\cd[0.62]{12}{1/3/B,2/5/B,4/9/B,6/10/B,7/11/B,8/12/B}\enskip
\cd[0.62]{12}{1/3/B,2/5/B,4/9/B,6/12/B,7/10/B,8/11/B}\\[4pt]
\cd[0.62]{12}{1/3/B,2/5/B,4/9/B,6/11/B,7/12/B,8/10/B}\enskip
\cd[0.62]{12}{1/3/B,2/5/B,4/11/B,6/9/B,7/10/B,8/12/B}\enskip
\cd[0.62]{12}{1/3/B,2/5/B,4/11/B,6/9/B,7/12/B,8/10/B}\enskip
\cd[0.62]{12}{1/3/B,2/5/B,4/10/B,6/11/B,7/9/B,8/12/B}\enskip
\cd[0.62]{12}{1/3/B,2/5/B,4/10/B,6/12/B,7/9/B,8/11/B}\enskip
\cd[0.62]{12}{1/3/B,2/6/B,4/7/B,5/10/B,8/11/B,9/12/B}\enskip
\cd[0.62]{12}{1/3/B,2/6/B,4/8/B,5/9/B,7/11/B,10/12/B}\enskip
\cd[0.62]{12}{1/3/B,2/6/B,4/8/B,5/11/B,7/9/B,10/12/B}\\[4pt]
\cd[0.62]{12}{1/3/B,2/6/B,4/8/B,5/10/B,7/12/B,9/11/B}\enskip
\cd[0.62]{12}{1/3/B,2/6/B,4/10/B,5/8/B,7/11/B,9/12/B}\enskip
\cd[0.62]{12}{1/3/B,2/6/B,4/10/B,5/8/B,7/12/B,9/11/B}\enskip
\cd[0.62]{12}{1/3/B,2/6/B,4/11/B,5/8/B,7/10/B,9/12/B}\enskip
\cd[0.62]{12}{1/3/B,2/6/B,4/9/B,5/10/B,7/11/B,8/12/B}\enskip
\cd[0.62]{12}{1/3/B,2/6/B,4/9/B,5/12/B,7/10/B,8/11/B}\enskip
\cd[0.62]{12}{1/3/B,2/6/B,4/9/B,5/11/B,7/12/B,8/10/B}\enskip
\cd[0.62]{12}{1/3/B,2/6/B,4/11/B,5/9/B,7/10/B,8/12/B}\\[4pt]
\cd[0.62]{12}{1/3/B,2/6/B,4/11/B,5/9/B,7/12/B,8/10/B}\enskip
\cd[0.62]{12}{1/3/B,2/6/B,4/10/B,5/11/B,7/9/B,8/12/B}\enskip
\cd[0.62]{12}{1/3/B,2/7/B,4/6/B,5/10/B,8/11/B,9/12/B}\enskip
\cd[0.62]{12}{1/3/B,2/8/B,4/6/B,5/10/B,7/11/B,9/12/B}\enskip
\cd[0.62]{12}{1/3/B,2/8/B,4/6/B,5/11/B,7/10/B,9/12/B}\enskip
\cd[0.62]{12}{1/3/B,2/9/B,4/6/B,5/10/B,7/11/B,8/12/B}\enskip
\cd[0.62]{12}{1/3/B,2/7/B,4/8/B,5/10/B,6/12/B,9/11/B}\enskip
\cd[0.62]{12}{1/3/B,2/7/B,4/10/B,5/8/B,6/11/B,9/12/B}\\[4pt]
\cd[0.62]{12}{1/3/B,2/7/B,4/10/B,5/8/B,6/12/B,9/11/B}\enskip
\cd[0.62]{12}{1/3/B,2/7/B,4/11/B,5/8/B,6/10/B,9/12/B}\enskip
\cd[0.62]{12}{1/3/B,2/7/B,4/9/B,5/10/B,6/11/B,8/12/B}\enskip
\cd[0.62]{12}{1/3/B,2/7/B,4/9/B,5/12/B,6/10/B,8/11/B}\enskip
\cd[0.62]{12}{1/3/B,2/7/B,4/9/B,5/11/B,6/12/B,8/10/B}\enskip
\cd[0.62]{12}{1/3/B,2/7/B,4/11/B,5/9/B,6/10/B,8/12/B}\enskip
\cd[0.62]{12}{1/3/B,2/7/B,4/10/B,5/11/B,6/9/B,8/12/B}\enskip
\cd[0.62]{12}{1/3/B,2/7/B,4/10/B,5/12/B,6/9/B,8/11/B}\\[4pt]
\cd[0.62]{12}{1/3/B,2/8/B,4/10/B,5/7/B,6/12/B,9/11/B}\enskip
\cd[0.62]{12}{1/3/B,2/8/B,4/9/B,5/10/B,6/11/B,7/12/B}\enskip
\cd[0.62]{12}{1/3/B,2/8/B,4/9/B,5/12/B,6/10/B,7/11/B}\enskip
\cd[0.62]{12}{1/3/B,2/8/B,4/9/B,5/11/B,6/12/B,7/10/B}\enskip
\cd[0.62]{12}{1/3/B,2/8/B,4/11/B,5/9/B,6/12/B,7/10/B}\enskip
\cd[0.62]{12}{1/3/B,2/8/B,4/11/B,5/12/B,6/9/B,7/10/B}\enskip
\cd[0.62]{12}{1/4/B,2/5/B,3/7/B,6/10/B,8/11/B,9/12/B}\enskip
\cd[0.62]{12}{1/4/B,2/5/B,3/9/B,6/10/B,7/11/B,8/12/B}\\[4pt]
\cd[0.62]{12}{1/4/B,2/5/B,3/9/B,6/12/B,7/10/B,8/11/B}\enskip
\cd[0.62]{12}{1/4/B,2/6/B,3/8/B,5/10/B,7/11/B,9/12/B}\enskip
\cd[0.62]{12}{1/4/B,2/6/B,3/8/B,5/11/B,7/10/B,9/12/B}\enskip
\cd[0.62]{12}{1/4/B,2/6/B,3/9/B,5/11/B,7/10/B,8/12/B}\enskip
\cd[0.62]{12}{1/4/B,2/6/B,3/10/B,5/9/B,7/11/B,8/12/B}\enskip
\cd[0.62]{12}{1/4/B,2/8/B,3/6/B,5/11/B,7/10/B,9/12/B}\enskip
\cd[0.62]{12}{1/4/B,2/10/B,3/6/B,5/9/B,7/11/B,8/12/B}\enskip
\cd[0.62]{12}{1/4/B,2/7/B,3/9/B,5/10/B,6/11/B,8/12/B}\\[4pt]
\cd[0.62]{12}{1/4/B,2/7/B,3/9/B,5/10/B,6/12/B,8/11/B}\enskip
\cd[0.62]{12}{1/4/B,2/7/B,3/10/B,5/9/B,6/12/B,8/11/B}\enskip
\cd[0.62]{12}{1/4/B,2/9/B,3/7/B,5/10/B,6/11/B,8/12/B}\enskip
\cd[0.62]{12}{1/4/B,2/10/B,3/7/B,5/9/B,6/12/B,8/11/B}\enskip
\cd[0.62]{12}{1/5/B,2/6/B,3/9/B,4/10/B,7/11/B,8/12/B}\enskip
\cd[0.62]{12}{1/5/B,2/7/B,3/9/B,4/10/B,6/11/B,8/12/B}\enskip
\cd[0.62]{12}{1/5/B,2/7/B,3/11/B,4/9/B,6/10/B,8/12/B}\enskip
\cd[0.62]{12}{1/5/B,2/8/B,3/9/B,4/10/B,6/11/B,7/12/B}\\[4pt]
\cd[0.62]{12}{1/5/B,2/10/B,3/8/B,4/9/B,6/11/B,7/12/B}\enskip
\cd[0.62]{12}{1/7/B,2/8/B,3/9/B,4/10/B,5/11/B,6/12/B}
\end{center}
\caption{One-cycle chord diagrams for genus $g=3$ ($k=6$), diagrams 1--80 of 82.}
\label{fig:atlas_g3_1}
\end{figure}

{\small
\paragraph{Duality, $g=3$.} 82 diagrams in total: 30 self-dual, 26 dual pairs.
\emph{Self-dual:} 1, 4, 5, 6, 9, 10, 13, 18, 20, 22, 24, 29, 31, 32, 37, 42, 46, 47, 50, 53, 56, 57, 58, 62, 65, 68, 70, 75, 77, 82.
\emph{Dual pairs:} (2,8), (3,11), (7,12), (14,41), (15,54), (16,35), (17,33), (19,44), (21,39), (23,51), (25,43), (26,55), (27,45), (28,49), (30,59), (34,52), (36,61), (38,60), (40,48), (63,81), (64,80), (66,79), (67,74), (69,72), (71,78), (73,76).
}

\subsection{Non-orientable Surfaces}

\begin{figure}[ht!]
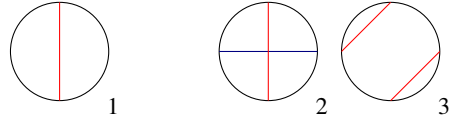

\begin{center}
\setcounter{cdn}{0}
\cd[0.65]{2}{1/2/R} \ \ \ \ \ \ \ \ \ \ \ \ \enskip
\cd[0.65]{4}{1/3/B,2/4/R}\enskip
\cd[0.65]{4}{1/2/R,3/4/R}
\end{center}
\caption{One-cycle chord diagrams for the projective plane ($g=-1$, $k=1$) and the Klein bottle ($g=-2$, $k=2$).}
\label{fig:atlas_gm1_1}
\end{figure}

{\small
\paragraph{Duality, $g=-1$.} 1 diagram in total: 1 self-dual, 0 dual pairs.
\emph{Self-dual:} 1.
}

{\small
\paragraph{Duality, $g=-2$.} 2 diagrams in total: 2 self-dual, 0 dual pairs.
\emph{Self-dual:} 2, 3.
}

\begin{figure}[ht!]
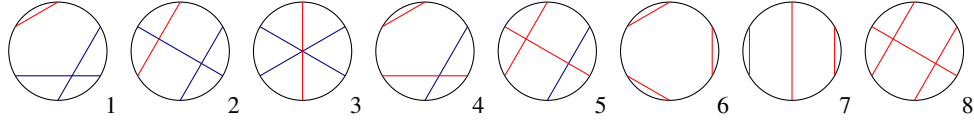

\begin{center}
\setcounter{cdn}{0}
\cd[0.65]{6}{1/3/B,2/4/B,5/6/R}\enskip
\cd[0.65]{6}{1/3/B,2/5/B,4/6/R}\enskip
\cd[0.65]{6}{1/4/B,2/5/B,3/6/R}\enskip
\cd[0.65]{6}{1/3/B,2/4/R,5/6/R}\enskip
\cd[0.65]{6}{1/3/B,2/5/R,4/6/R}\enskip
\cd[0.65]{6}{1/2/R,3/4/R,5/6/R}\enskip
\cd[0.65]{6}{1/2/R,3/6/R,4/5/R}\enskip
\cd[0.65]{6}{1/3/R,2/5/R,4/6/R}
\end{center}
\caption{One-cycle chord diagrams for non-orientable genus $g=-3$ ($k=3$).}
\label{fig:atlas_gm3_1}
\end{figure}

{\small
\paragraph{Duality, $g=-3$.} 8 diagrams in total: 4 self-dual, 2 dual pairs.
\emph{Self-dual:} 1, 4, 6, 7.
\emph{Dual pairs:} (2,5), (3,8).
}
\begin{figure}[ht!]
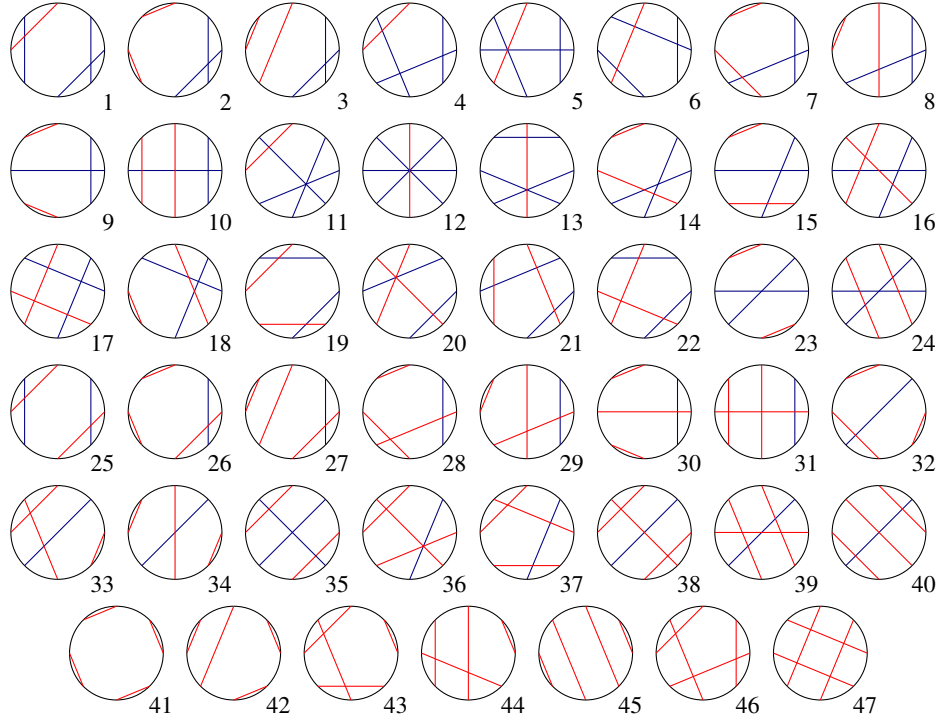

\begin{center}
\setcounter{cdn}{0}
\cd[0.62]{8}{1/3/B,2/4/B,5/7/B,6/8/R}\enskip
\cd[0.62]{8}{1/3/B,2/4/B,5/6/R,7/8/R}\enskip
\cd[0.62]{8}{1/3/B,2/4/B,5/8/R,6/7/R}\enskip
\cd[0.62]{8}{1/3/B,2/5/B,4/7/B,6/8/R}\enskip
\cd[0.62]{8}{1/3/B,2/6/B,4/7/B,5/8/R}\enskip
\cd[0.62]{8}{1/3/B,2/7/B,4/6/B,5/8/R}\enskip
\cd[0.62]{8}{1/3/B,2/5/B,4/6/R,7/8/R}\enskip
\cd[0.62]{8}{1/3/B,2/5/B,4/8/R,6/7/R}\\[4pt]
\cd[0.62]{8}{1/3/B,2/6/B,4/5/R,7/8/R}\enskip
\cd[0.62]{8}{1/3/B,2/6/B,4/8/R,5/7/R}\enskip
\cd[0.62]{8}{1/4/B,2/5/B,3/7/B,6/8/R}\enskip
\cd[0.62]{8}{1/5/B,2/6/B,3/7/B,4/8/R}\enskip
\cd[0.62]{8}{1/7/B,2/5/B,3/6/B,4/8/R}\enskip
\cd[0.62]{8}{1/4/B,2/5/B,3/6/R,7/8/R}\enskip
\cd[0.62]{8}{1/4/B,2/6/B,3/5/R,7/8/R}\enskip
\cd[0.62]{8}{1/4/B,2/6/B,3/7/R,5/8/R}\\[4pt]
\cd[0.62]{8}{1/4/B,2/7/B,3/6/R,5/8/R}\enskip
\cd[0.62]{8}{1/4/B,2/7/B,3/8/R,5/6/R}\enskip
\cd[0.62]{8}{1/7/B,2/4/B,3/5/R,6/8/R}\enskip
\cd[0.62]{8}{1/6/B,2/4/B,3/7/R,5/8/R}\enskip
\cd[0.62]{8}{1/6/B,2/4/B,3/8/R,5/7/R}\enskip
\cd[0.62]{8}{1/7/B,2/4/B,3/6/R,5/8/R}\enskip
\cd[0.62]{8}{1/5/B,2/6/B,3/4/R,7/8/R}\enskip
\cd[0.62]{8}{1/5/B,2/6/B,3/8/R,4/7/R}\\[4pt]
\cd[0.62]{8}{1/3/B,2/4/R,5/7/B,6/8/R}\enskip
\cd[0.62]{8}{1/3/B,2/4/R,5/6/R,7/8/R}\enskip
\cd[0.62]{8}{1/3/B,2/4/R,5/8/R,6/7/R}\enskip
\cd[0.62]{8}{1/3/B,2/5/R,4/6/R,7/8/R}\enskip
\cd[0.62]{8}{1/3/B,2/5/R,4/8/R,6/7/R}\enskip
\cd[0.62]{8}{1/3/B,2/6/R,4/5/R,7/8/R}\enskip
\cd[0.62]{8}{1/3/B,2/6/R,4/8/R,5/7/R}\enskip
\cd[0.62]{8}{1/5/B,2/3/R,4/6/R,7/8/R}\\[4pt]
\cd[0.62]{8}{1/5/B,2/3/R,4/7/R,6/8/R}\enskip
\cd[0.62]{8}{1/5/B,2/3/R,4/8/R,6/7/R}\enskip
\cd[0.62]{8}{1/5/B,2/4/R,3/7/B,6/8/R}\enskip
\cd[0.62]{8}{1/4/B,2/5/R,3/7/R,6/8/R}\enskip
\cd[0.62]{8}{1/4/B,2/7/R,3/5/R,6/8/R}\enskip
\cd[0.62]{8}{1/5/B,2/4/R,3/7/R,6/8/R}\enskip
\cd[0.62]{8}{1/5/B,2/6/R,3/8/R,4/7/R}\enskip
\cd[0.62]{8}{1/5/B,2/8/R,3/7/R,4/6/R}\\[4pt]
\cd[0.62]{8}{1/2/R,3/4/R,5/6/R,7/8/R}\enskip
\cd[0.62]{8}{1/2/R,3/4/R,5/8/R,6/7/R}\enskip
\cd[0.62]{8}{1/2/R,3/5/R,4/7/R,6/8/R}\enskip
\cd[0.62]{8}{1/2/R,3/6/R,4/8/R,5/7/R}\enskip
\cd[0.62]{8}{1/2/R,3/8/R,4/7/R,5/6/R}\enskip
\cd[0.62]{8}{1/3/R,2/5/R,4/7/R,6/8/R}\enskip
\cd[0.62]{8}{1/4/R,2/7/R,3/6/R,5/8/R}
\end{center}
\caption{One-cycle chord diagrams for non-orientable genus $g=-4$ ($k=4$).}
\label{fig:atlas_gm4_1}
\end{figure}

{\small
\paragraph{Duality, $g=-4$.} 47 diagrams in total: 23 self-dual, 12 dual pairs.
\emph{Self-dual:} 1, 2, 3, 6, 9, 10, 13, 19, 22, 23, 24, 25, 26, 27, 30, 31, 35, 36, 38, 41, 42, 45, 47.
\emph{Dual pairs:} (4,20), (5,21), (7,29), (8,28), (11,39), (12,40), (14,44), (15,33), (16,37), (17,46), (18,43), (32,34).
}

\clearpage
\begin{figure}[ht!]
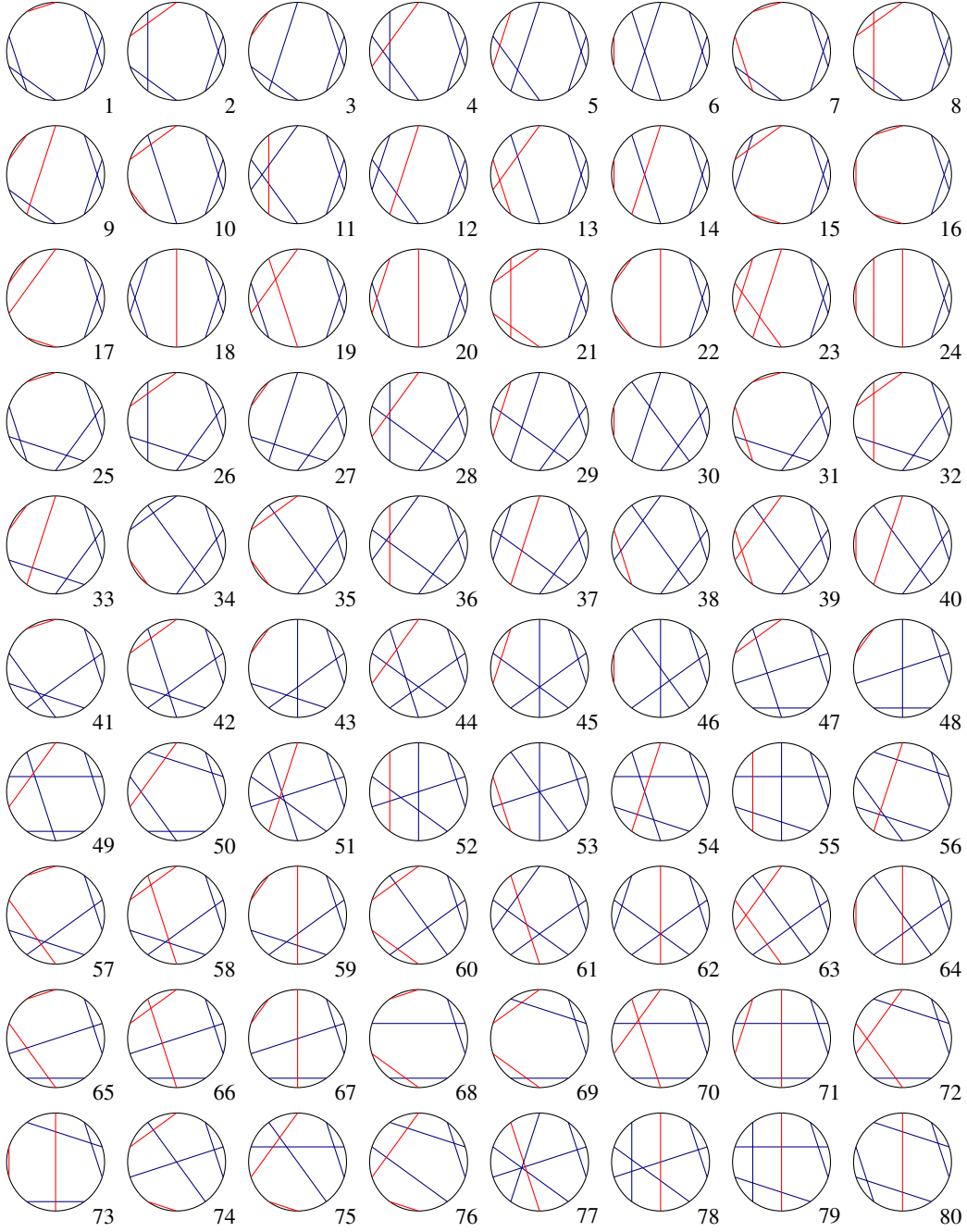

\begin{center}
\setcounter{cdn}{0}
\cd[0.7]{10}{1/3/B,2/4/B,5/7/B,6/8/B,9/10/R}\enskip
\cd[0.7]{10}{1/3/B,2/4/B,5/7/B,6/9/B,8/10/R}\enskip
\cd[0.7]{10}{1/3/B,2/4/B,5/7/B,6/10/B,8/9/R}\enskip
\cd[0.7]{10}{1/3/B,2/4/B,5/8/B,6/9/B,7/10/R}\enskip
\cd[0.7]{10}{1/3/B,2/4/B,5/8/B,6/10/B,7/9/R}\enskip
\cd[0.7]{10}{1/3/B,2/4/B,5/9/B,6/10/B,7/8/R}\enskip
\cd[0.7]{10}{1/3/B,2/4/B,5/7/B,6/8/R,9/10/R}\enskip
\cd[0.7]{10}{1/3/B,2/4/B,5/7/B,6/9/R,8/10/R}\\[4pt]
\cd[0.7]{10}{1/3/B,2/4/B,5/7/B,6/10/R,8/9/R}\enskip
\cd[0.7]{10}{1/3/B,2/4/B,5/9/B,6/7/R,8/10/R}\enskip
\cd[0.7]{10}{1/3/B,2/4/B,5/8/B,6/9/R,7/10/B}\enskip
\cd[0.7]{10}{1/3/B,2/4/B,5/8/B,6/10/R,7/9/B}\enskip
\cd[0.7]{10}{1/3/B,2/4/B,5/9/B,6/8/R,7/10/R}\enskip
\cd[0.7]{10}{1/3/B,2/4/B,5/9/B,6/10/R,7/8/R}\enskip
\cd[0.7]{10}{1/3/B,2/4/B,5/6/R,7/9/B,8/10/R}\enskip
\cd[0.7]{10}{1/3/B,2/4/B,5/6/R,7/8/R,9/10/R}\\[4pt]
\cd[0.7]{10}{1/3/B,2/4/B,5/6/R,7/10/R,8/9/R}\enskip
\cd[0.7]{10}{1/3/B,2/4/B,5/10/R,6/8/B,7/9/B}\enskip
\cd[0.7]{10}{1/3/B,2/4/B,5/9/R,6/8/B,7/10/R}\enskip
\cd[0.7]{10}{1/3/B,2/4/B,5/10/R,6/8/B,7/9/R}\enskip
\cd[0.7]{10}{1/3/B,2/4/B,5/7/R,6/9/R,8/10/R}\enskip
\cd[0.7]{10}{1/3/B,2/4/B,5/10/R,6/7/R,8/9/R}\enskip
\cd[0.7]{10}{1/3/B,2/4/B,5/8/R,6/10/R,7/9/R}\enskip
\cd[0.7]{10}{1/3/B,2/4/B,5/10/R,6/9/R,7/8/R}\\[4pt]
\cd[0.7]{10}{1/3/B,2/5/B,4/7/B,6/8/B,9/10/R}\enskip
\cd[0.7]{10}{1/3/B,2/5/B,4/7/B,6/9/B,8/10/R}\enskip
\cd[0.7]{10}{1/3/B,2/5/B,4/7/B,6/10/B,8/9/R}\enskip
\cd[0.7]{10}{1/3/B,2/5/B,4/8/B,6/9/B,7/10/R}\enskip
\cd[0.7]{10}{1/3/B,2/5/B,4/8/B,6/10/B,7/9/R}\enskip
\cd[0.7]{10}{1/3/B,2/5/B,4/9/B,6/10/B,7/8/R}\enskip
\cd[0.7]{10}{1/3/B,2/5/B,4/7/B,6/8/R,9/10/R}\enskip
\cd[0.7]{10}{1/3/B,2/5/B,4/7/B,6/9/R,8/10/R}\\[4pt]
\cd[0.7]{10}{1/3/B,2/5/B,4/7/B,6/10/R,8/9/R}\enskip
\cd[0.7]{10}{1/3/B,2/5/B,4/9/B,6/7/R,8/10/B}\enskip
\cd[0.7]{10}{1/3/B,2/5/B,4/9/B,6/7/R,8/10/R}\enskip
\cd[0.7]{10}{1/3/B,2/5/B,4/8/B,6/9/R,7/10/B}\enskip
\cd[0.7]{10}{1/3/B,2/5/B,4/8/B,6/10/R,7/9/B}\enskip
\cd[0.7]{10}{1/3/B,2/5/B,4/9/B,6/8/R,7/10/B}\enskip
\cd[0.7]{10}{1/3/B,2/5/B,4/9/B,6/8/R,7/10/R}\enskip
\cd[0.7]{10}{1/3/B,2/5/B,4/9/B,6/10/R,7/8/R}\\[4pt]
\cd[0.7]{10}{1/3/B,2/6/B,4/7/B,5/8/B,9/10/R}\enskip
\cd[0.7]{10}{1/3/B,2/6/B,4/7/B,5/9/B,8/10/R}\enskip
\cd[0.7]{10}{1/3/B,2/6/B,4/7/B,5/10/B,8/9/R}\enskip
\cd[0.7]{10}{1/3/B,2/6/B,4/8/B,5/9/B,7/10/R}\enskip
\cd[0.7]{10}{1/3/B,2/6/B,4/8/B,5/10/B,7/9/R}\enskip
\cd[0.7]{10}{1/3/B,2/6/B,4/9/B,5/10/B,7/8/R}\enskip
\cd[0.7]{10}{1/3/B,2/7/B,4/6/B,5/9/B,8/10/R}\enskip
\cd[0.7]{10}{1/3/B,2/7/B,4/6/B,5/10/B,8/9/R}\\[4pt]
\cd[0.7]{10}{1/3/B,2/8/B,4/6/B,5/9/B,7/10/R}\enskip
\cd[0.7]{10}{1/3/B,2/9/B,4/6/B,5/8/B,7/10/R}\enskip
\cd[0.7]{10}{1/3/B,2/7/B,4/8/B,5/9/B,6/10/R}\enskip
\cd[0.7]{10}{1/3/B,2/7/B,4/8/B,5/10/B,6/9/R}\enskip
\cd[0.7]{10}{1/3/B,2/7/B,4/9/B,5/10/B,6/8/R}\enskip
\cd[0.7]{10}{1/3/B,2/8/B,4/7/B,5/9/B,6/10/R}\enskip
\cd[0.7]{10}{1/3/B,2/8/B,4/7/B,5/10/B,6/9/R}\enskip
\cd[0.7]{10}{1/3/B,2/9/B,4/7/B,5/8/B,6/10/R}\\[4pt]
\cd[0.7]{10}{1/3/B,2/6/B,4/7/B,5/8/R,9/10/R}\enskip
\cd[0.7]{10}{1/3/B,2/6/B,4/7/B,5/9/R,8/10/R}\enskip
\cd[0.7]{10}{1/3/B,2/6/B,4/7/B,5/10/R,8/9/R}\enskip
\cd[0.7]{10}{1/3/B,2/6/B,4/9/B,5/7/R,8/10/R}\enskip
\cd[0.7]{10}{1/3/B,2/6/B,4/8/B,5/9/R,7/10/B}\enskip
\cd[0.7]{10}{1/3/B,2/6/B,4/8/B,5/10/R,7/9/B}\enskip
\cd[0.7]{10}{1/3/B,2/6/B,4/9/B,5/8/R,7/10/R}\enskip
\cd[0.7]{10}{1/3/B,2/6/B,4/9/B,5/10/R,7/8/R}\\[4pt]
\cd[0.7]{10}{1/3/B,2/7/B,4/6/B,5/8/R,9/10/R}\enskip
\cd[0.7]{10}{1/3/B,2/7/B,4/6/B,5/9/R,8/10/R}\enskip
\cd[0.7]{10}{1/3/B,2/7/B,4/6/B,5/10/R,8/9/R}\enskip
\cd[0.7]{10}{1/3/B,2/8/B,4/6/B,5/7/R,9/10/R}\enskip
\cd[0.7]{10}{1/3/B,2/9/B,4/6/B,5/7/R,8/10/R}\enskip
\cd[0.7]{10}{1/3/B,2/8/B,4/6/B,5/9/R,7/10/R}\enskip
\cd[0.7]{10}{1/3/B,2/8/B,4/6/B,5/10/R,7/9/R}\enskip
\cd[0.7]{10}{1/3/B,2/9/B,4/6/B,5/8/R,7/10/R}\\[4pt]
\cd[0.7]{10}{1/3/B,2/9/B,4/6/B,5/10/R,7/8/R}\enskip
\cd[0.7]{10}{1/3/B,2/7/B,4/9/B,5/6/R,8/10/R}\enskip
\cd[0.7]{10}{1/3/B,2/8/B,4/9/B,5/6/R,7/10/R}\enskip
\cd[0.7]{10}{1/3/B,2/9/B,4/8/B,5/6/R,7/10/R}\enskip
\cd[0.7]{10}{1/3/B,2/7/B,4/8/B,5/9/R,6/10/B}\enskip
\cd[0.7]{10}{1/3/B,2/7/B,4/8/B,5/10/R,6/9/B}\enskip
\cd[0.7]{10}{1/3/B,2/8/B,4/7/B,5/10/R,6/9/B}\enskip
\cd[0.7]{10}{1/3/B,2/9/B,4/7/B,5/10/R,6/8/B}
\end{center}
\caption{One-cycle chord diagrams for non-orientable genus $g=-5$ ($k=5$), diagrams 1--80 of 473.}
\label{fig:atlas_gm5_1}
\end{figure}

\begin{figure}[ht!]
\begin{center}
\cd[0.7]{10}{1/3/B,2/7/B,4/9/B,5/8/R,6/10/R}\enskip
\cd[0.7]{10}{1/3/B,2/7/B,4/9/B,5/10/R,6/8/R}\enskip
\cd[0.7]{10}{1/3/B,2/8/B,4/9/B,5/7/R,6/10/R}\enskip
\cd[0.7]{10}{1/3/B,2/9/B,4/8/B,5/7/R,6/10/R}\enskip
\cd[0.7]{10}{1/3/B,2/8/B,4/9/B,5/10/R,6/7/R}\enskip
\cd[0.7]{10}{1/3/B,2/9/B,4/8/B,5/10/R,6/7/R}\enskip
\cd[0.7]{10}{1/3/B,2/5/B,4/6/R,7/9/B,8/10/R}\enskip
\cd[0.7]{10}{1/3/B,2/5/B,4/6/R,7/8/R,9/10/R}\\[4pt]
\cd[0.7]{10}{1/3/B,2/5/B,4/6/R,7/10/R,8/9/R}\enskip
\cd[0.7]{10}{1/3/B,2/5/B,4/7/R,6/8/B,9/10/R}\enskip
\cd[0.7]{10}{1/3/B,2/5/B,4/9/R,6/8/B,7/10/R}\enskip
\cd[0.7]{10}{1/3/B,2/5/B,4/10/R,6/8/B,7/9/R}\enskip
\cd[0.7]{10}{1/3/B,2/5/B,4/7/R,6/9/R,8/10/R}\enskip
\cd[0.7]{10}{1/3/B,2/5/B,4/8/R,6/7/R,9/10/R}\enskip
\cd[0.7]{10}{1/3/B,2/5/B,4/10/R,6/7/R,8/9/R}\enskip
\cd[0.7]{10}{1/3/B,2/5/B,4/8/R,6/10/R,7/9/B}\\[4pt]
\cd[0.7]{10}{1/3/B,2/5/B,4/10/R,6/8/R,7/9/B}\enskip
\cd[0.7]{10}{1/3/B,2/5/B,4/8/R,6/10/R,7/9/R}\enskip
\cd[0.7]{10}{1/3/B,2/5/B,4/9/R,6/8/R,7/10/R}\enskip
\cd[0.7]{10}{1/3/B,2/5/B,4/10/R,6/9/R,7/8/R}\enskip
\cd[0.7]{10}{1/3/B,2/6/B,4/5/R,7/9/B,8/10/R}\enskip
\cd[0.7]{10}{1/3/B,2/6/B,4/5/R,7/8/R,9/10/R}\enskip
\cd[0.7]{10}{1/3/B,2/6/B,4/5/R,7/10/R,8/9/R}\enskip
\cd[0.7]{10}{1/3/B,2/7/B,4/5/R,6/9/B,8/10/R}\\[4pt]
\cd[0.7]{10}{1/3/B,2/8/B,4/5/R,6/9/B,7/10/R}\enskip
\cd[0.7]{10}{1/3/B,2/7/B,4/5/R,6/8/R,9/10/R}\enskip
\cd[0.7]{10}{1/3/B,2/7/B,4/5/R,6/10/R,8/9/R}\enskip
\cd[0.7]{10}{1/3/B,2/8/B,4/5/R,6/10/R,7/9/R}\enskip
\cd[0.7]{10}{1/3/B,2/6/B,4/9/R,5/8/B,7/10/R}\enskip
\cd[0.7]{10}{1/3/B,2/6/B,4/10/R,5/8/B,7/9/R}\enskip
\cd[0.7]{10}{1/3/B,2/7/B,4/10/R,5/8/B,6/9/B}\enskip
\cd[0.7]{10}{1/3/B,2/7/B,4/9/R,5/8/B,6/10/R}\\[4pt]
\cd[0.7]{10}{1/3/B,2/7/B,4/10/R,5/8/B,6/9/R}\enskip
\cd[0.7]{10}{1/3/B,2/8/B,4/9/R,5/7/B,6/10/R}\enskip
\cd[0.7]{10}{1/3/B,2/8/B,4/10/R,5/7/B,6/9/R}\enskip
\cd[0.7]{10}{1/3/B,2/6/B,4/7/R,5/9/R,8/10/R}\enskip
\cd[0.7]{10}{1/3/B,2/6/B,4/10/R,5/7/R,8/9/R}\enskip
\cd[0.7]{10}{1/3/B,2/6/B,4/8/R,5/10/R,7/9/R}\enskip
\cd[0.7]{10}{1/3/B,2/6/B,4/9/R,5/8/R,7/10/R}\enskip
\cd[0.7]{10}{1/3/B,2/6/B,4/10/R,5/9/R,7/8/R}\\[4pt]
\cd[0.7]{10}{1/3/B,2/7/B,4/6/R,5/9/R,8/10/R}\enskip
\cd[0.7]{10}{1/3/B,2/7/B,4/10/R,5/6/R,8/9/R}\enskip
\cd[0.7]{10}{1/3/B,2/7/B,4/8/R,5/10/R,6/9/R}\enskip
\cd[0.7]{10}{1/3/B,2/7/B,4/10/R,5/9/R,6/8/R}\enskip
\cd[0.7]{10}{1/4/B,2/5/B,3/7/B,6/9/B,8/10/R}\enskip
\cd[0.7]{10}{1/4/B,2/5/B,3/7/B,6/10/B,8/9/R}\enskip
\cd[0.7]{10}{1/4/B,2/5/B,3/8/B,6/9/B,7/10/R}\enskip
\cd[0.7]{10}{1/4/B,2/5/B,3/8/B,6/10/B,7/9/R}\\[4pt]
\cd[0.7]{10}{1/4/B,2/5/B,3/7/B,6/8/R,9/10/R}\enskip
\cd[0.7]{10}{1/4/B,2/5/B,3/7/B,6/9/R,8/10/R}\enskip
\cd[0.7]{10}{1/4/B,2/5/B,3/7/B,6/10/R,8/9/R}\enskip
\cd[0.7]{10}{1/4/B,2/6/B,3/9/B,5/8/B,7/10/R}\enskip
\cd[0.7]{10}{1/4/B,2/9/B,3/6/B,5/8/B,7/10/R}\enskip
\cd[0.7]{10}{1/4/B,2/7/B,3/9/B,5/8/B,6/10/R}\enskip
\cd[0.7]{10}{1/4/B,2/9/B,3/7/B,5/8/B,6/10/R}\enskip
\cd[0.7]{10}{1/4/B,2/6/B,3/8/B,5/7/R,9/10/R}\\[4pt]
\cd[0.7]{10}{1/4/B,2/6/B,3/8/B,5/9/R,7/10/R}\enskip
\cd[0.7]{10}{1/4/B,2/6/B,3/8/B,5/10/R,7/9/R}\enskip
\cd[0.7]{10}{1/4/B,2/8/B,3/6/B,5/7/R,9/10/R}\enskip
\cd[0.7]{10}{1/4/B,2/8/B,3/6/B,5/9/R,7/10/R}\enskip
\cd[0.7]{10}{1/4/B,2/8/B,3/6/B,5/10/R,7/9/R}\enskip
\cd[0.7]{10}{1/5/B,2/6/B,3/7/B,4/8/B,9/10/R}\enskip
\cd[0.7]{10}{1/5/B,2/6/B,3/7/B,4/9/B,8/10/R}\enskip
\cd[0.7]{10}{1/5/B,2/6/B,3/8/B,4/9/B,7/10/R}\\[4pt]
\cd[0.7]{10}{1/5/B,2/7/B,3/8/B,4/9/B,6/10/R}\enskip
\cd[0.7]{10}{1/5/B,2/9/B,3/7/B,4/8/B,6/10/R}\enskip
\cd[0.7]{10}{1/6/B,2/7/B,3/8/B,4/9/B,5/10/R}\enskip
\cd[0.7]{10}{1/6/B,2/9/B,3/7/B,4/8/B,5/10/R}\enskip
\cd[0.7]{10}{1/8/B,2/6/B,3/9/B,4/7/B,5/10/R}\enskip
\cd[0.7]{10}{1/8/B,2/9/B,3/6/B,4/7/B,5/10/R}\enskip
\cd[0.7]{10}{1/5/B,2/6/B,3/7/B,4/8/R,9/10/R}\enskip
\cd[0.7]{10}{1/5/B,2/6/B,3/7/B,4/9/R,8/10/R}\\[4pt]
\cd[0.7]{10}{1/5/B,2/6/B,3/9/B,4/7/R,8/10/R}\enskip
\cd[0.7]{10}{1/5/B,2/6/B,3/9/B,4/8/R,7/10/R}\enskip
\cd[0.7]{10}{1/5/B,2/6/B,3/9/B,4/10/R,7/8/R}\enskip
\cd[0.7]{10}{1/5/B,2/8/B,3/6/B,4/7/R,9/10/R}\enskip
\cd[0.7]{10}{1/5/B,2/8/B,3/6/B,4/9/R,7/10/R}\enskip
\cd[0.7]{10}{1/5/B,2/8/B,3/6/B,4/10/R,7/9/R}\enskip
\cd[0.7]{10}{1/5/B,2/7/B,3/8/B,4/6/R,9/10/R}\enskip
\cd[0.7]{10}{1/5/B,2/9/B,3/7/B,4/6/R,8/10/R}
\end{center}
\caption{One-cycle chord diagrams for non-orientable genus $g=-5$ ($k=5$), diagrams 81--160 of 473.}
\label{fig:atlas_gm5_2}
\end{figure}

\begin{figure}[ht!]
\begin{center}
\cd[0.7]{10}{1/5/B,2/7/B,3/8/B,4/9/R,6/10/R}\enskip
\cd[0.7]{10}{1/5/B,2/7/B,3/8/B,4/10/R,6/9/R}\enskip
\cd[0.7]{10}{1/5/B,2/9/B,3/7/B,4/8/R,6/10/R}\enskip
\cd[0.7]{10}{1/5/B,2/9/B,3/7/B,4/10/R,6/8/R}\enskip
\cd[0.7]{10}{1/7/B,2/5/B,3/6/B,4/8/R,9/10/R}\enskip
\cd[0.7]{10}{1/7/B,2/5/B,3/6/B,4/9/R,8/10/R}\enskip
\cd[0.7]{10}{1/9/B,2/5/B,3/6/B,4/7/R,8/10/R}\enskip
\cd[0.7]{10}{1/9/B,2/5/B,3/6/B,4/8/R,7/10/R}\\[4pt]
\cd[0.7]{10}{1/9/B,2/5/B,3/6/B,4/10/R,7/8/R}\enskip
\cd[0.7]{10}{1/7/B,2/5/B,3/8/B,4/6/R,9/10/R}\enskip
\cd[0.7]{10}{1/9/B,2/5/B,3/7/B,4/6/R,8/10/R}\enskip
\cd[0.7]{10}{1/9/B,2/5/B,3/8/B,4/6/R,7/10/R}\enskip
\cd[0.7]{10}{1/7/B,2/5/B,3/8/B,4/9/R,6/10/R}\enskip
\cd[0.7]{10}{1/7/B,2/5/B,3/8/B,4/10/R,6/9/R}\enskip
\cd[0.7]{10}{1/9/B,2/5/B,3/7/B,4/8/R,6/10/R}\enskip
\cd[0.7]{10}{1/9/B,2/5/B,3/7/B,4/10/R,6/8/R}\\[4pt]
\cd[0.7]{10}{1/9/B,2/5/B,3/8/B,4/7/R,6/10/R}\enskip
\cd[0.7]{10}{1/9/B,2/5/B,3/8/B,4/10/R,6/7/R}\enskip
\cd[0.7]{10}{1/6/B,2/7/B,3/5/B,4/8/R,9/10/R}\enskip
\cd[0.7]{10}{1/6/B,2/7/B,3/5/B,4/9/R,8/10/R}\enskip
\cd[0.7]{10}{1/6/B,2/8/B,3/5/B,4/7/R,9/10/R}\enskip
\cd[0.7]{10}{1/6/B,2/8/B,3/5/B,4/9/R,7/10/R}\enskip
\cd[0.7]{10}{1/6/B,2/8/B,3/5/B,4/10/R,7/9/R}\enskip
\cd[0.7]{10}{1/7/B,2/8/B,3/5/B,4/6/R,9/10/R}\\[4pt]
\cd[0.7]{10}{1/7/B,2/8/B,3/5/B,4/9/R,6/10/R}\enskip
\cd[0.7]{10}{1/7/B,2/8/B,3/5/B,4/10/R,6/9/R}\enskip
\cd[0.7]{10}{1/4/B,2/5/B,3/6/R,7/9/B,8/10/R}\enskip
\cd[0.7]{10}{1/4/B,2/5/B,3/6/R,7/8/R,9/10/R}\enskip
\cd[0.7]{10}{1/4/B,2/5/B,3/6/R,7/10/R,8/9/R}\enskip
\cd[0.7]{10}{1/4/B,2/5/B,3/7/R,6/9/R,8/10/R}\enskip
\cd[0.7]{10}{1/4/B,2/5/B,3/8/R,6/7/R,9/10/R}\enskip
\cd[0.7]{10}{1/4/B,2/5/B,3/8/R,6/10/R,7/9/B}\\[4pt]
\cd[0.7]{10}{1/4/B,2/5/B,3/8/R,6/10/R,7/9/R}\enskip
\cd[0.7]{10}{1/4/B,2/6/B,3/5/R,7/9/B,8/10/R}\enskip
\cd[0.7]{10}{1/4/B,2/6/B,3/5/R,7/8/R,9/10/R}\enskip
\cd[0.7]{10}{1/4/B,2/6/B,3/5/R,7/10/R,8/9/R}\enskip
\cd[0.7]{10}{1/4/B,2/9/B,3/5/R,6/8/B,7/10/R}\enskip
\cd[0.7]{10}{1/4/B,2/7/B,3/5/R,6/9/R,8/10/R}\enskip
\cd[0.7]{10}{1/4/B,2/8/B,3/5/R,6/7/R,9/10/R}\enskip
\cd[0.7]{10}{1/4/B,2/8/B,3/5/R,6/10/R,7/9/R}\\[4pt]
\cd[0.7]{10}{1/4/B,2/9/B,3/5/R,6/8/R,7/10/R}\enskip
\cd[0.7]{10}{1/4/B,2/7/B,3/6/R,5/9/B,8/10/R}\enskip
\cd[0.7]{10}{1/4/B,2/8/B,3/6/R,5/7/B,9/10/R}\enskip
\cd[0.7]{10}{1/4/B,2/8/B,3/6/R,5/9/B,7/10/R}\enskip
\cd[0.7]{10}{1/4/B,2/7/B,3/8/R,5/9/B,6/10/R}\enskip
\cd[0.7]{10}{1/4/B,2/7/B,3/10/R,5/9/B,6/8/R}\enskip
\cd[0.7]{10}{1/4/B,2/8/B,3/7/R,5/9/B,6/10/R}\enskip
\cd[0.7]{10}{1/4/B,2/8/B,3/9/R,5/7/B,6/10/R}\\[4pt]
\cd[0.7]{10}{1/4/B,2/8/B,3/10/R,5/7/B,6/9/R}\enskip
\cd[0.7]{10}{1/4/B,2/9/B,3/8/R,5/7/B,6/10/R}\enskip
\cd[0.7]{10}{1/4/B,2/9/B,3/10/R,5/7/B,6/8/R}\enskip
\cd[0.7]{10}{1/4/B,2/8/B,3/10/R,5/9/B,6/7/R}\enskip
\cd[0.7]{10}{1/4/B,2/6/B,3/7/R,5/8/R,9/10/R}\enskip
\cd[0.7]{10}{1/4/B,2/6/B,3/7/R,5/10/R,8/9/R}\enskip
\cd[0.7]{10}{1/4/B,2/6/B,3/9/R,5/7/R,8/10/R}\enskip
\cd[0.7]{10}{1/4/B,2/6/B,3/8/R,5/10/R,7/9/B}\\[4pt]
\cd[0.7]{10}{1/4/B,2/6/B,3/10/R,5/8/R,7/9/B}\enskip
\cd[0.7]{10}{1/4/B,2/6/B,3/8/R,5/9/R,7/10/R}\enskip
\cd[0.7]{10}{1/4/B,2/6/B,3/10/R,5/8/R,7/9/R}\enskip
\cd[0.7]{10}{1/4/B,2/6/B,3/9/R,5/10/R,7/8/R}\enskip
\cd[0.7]{10}{1/4/B,2/7/B,3/6/R,5/8/R,9/10/R}\enskip
\cd[0.7]{10}{1/4/B,2/7/B,3/6/R,5/10/R,8/9/R}\enskip
\cd[0.7]{10}{1/4/B,2/8/B,3/6/R,5/9/R,7/10/R}\enskip
\cd[0.7]{10}{1/4/B,2/9/B,3/6/R,5/10/R,7/8/R}\\[4pt]
\cd[0.7]{10}{1/4/B,2/7/B,3/8/R,5/6/R,9/10/R}\enskip
\cd[0.7]{10}{1/4/B,2/7/B,3/10/R,5/6/R,8/9/R}\enskip
\cd[0.7]{10}{1/4/B,2/8/B,3/7/R,5/6/R,9/10/R}\enskip
\cd[0.7]{10}{1/4/B,2/8/B,3/9/R,5/6/R,7/10/R}\enskip
\cd[0.7]{10}{1/4/B,2/9/B,3/10/R,5/6/R,7/8/R}\enskip
\cd[0.7]{10}{1/4/B,2/7/B,3/8/R,5/10/R,6/9/R}\enskip
\cd[0.7]{10}{1/4/B,2/7/B,3/9/R,5/8/R,6/10/R}\enskip
\cd[0.7]{10}{1/4/B,2/7/B,3/10/R,5/9/R,6/8/R}\\[4pt]
\cd[0.7]{10}{1/4/B,2/8/B,3/7/R,5/10/R,6/9/R}\enskip
\cd[0.7]{10}{1/4/B,2/9/B,3/7/R,5/8/R,6/10/R}\enskip
\cd[0.7]{10}{1/4/B,2/8/B,3/10/R,5/7/R,6/9/R}\enskip
\cd[0.7]{10}{1/4/B,2/8/B,3/9/R,5/10/R,6/7/R}\enskip
\cd[0.7]{10}{1/4/B,2/9/B,3/8/R,5/10/R,6/7/R}\enskip
\cd[0.7]{10}{1/4/B,2/9/B,3/10/R,5/8/R,6/7/R}\enskip
\cd[0.7]{10}{1/7/B,2/4/B,3/5/R,6/8/R,9/10/R}\enskip
\cd[0.7]{10}{1/7/B,2/4/B,3/5/R,6/9/R,8/10/R}
\end{center}
\caption{One-cycle chord diagrams for non-orientable genus $g=-5$ ($k=5$), diagrams 161--240 of 473.}
\label{fig:atlas_gm5_3}
\end{figure}

\begin{figure}[ht!]
\begin{center}
\cd[0.7]{10}{1/7/B,2/4/B,3/5/R,6/10/R,8/9/R}\enskip
\cd[0.7]{10}{1/9/B,2/4/B,3/5/R,6/7/R,8/10/R}\enskip
\cd[0.7]{10}{1/9/B,2/4/B,3/5/R,6/8/R,7/10/R}\enskip
\cd[0.7]{10}{1/9/B,2/4/B,3/5/R,6/10/R,7/8/R}\enskip
\cd[0.7]{10}{1/7/B,2/4/B,3/6/R,5/9/B,8/10/R}\enskip
\cd[0.7]{10}{1/7/B,2/4/B,3/8/R,5/9/B,6/10/R}\enskip
\cd[0.7]{10}{1/6/B,2/4/B,3/7/R,5/8/R,9/10/R}\enskip
\cd[0.7]{10}{1/6/B,2/4/B,3/7/R,5/9/R,8/10/R}\\[4pt]
\cd[0.7]{10}{1/6/B,2/4/B,3/7/R,5/10/R,8/9/R}\enskip
\cd[0.7]{10}{1/6/B,2/4/B,3/8/R,5/7/R,9/10/R}\enskip
\cd[0.7]{10}{1/6/B,2/4/B,3/9/R,5/7/R,8/10/R}\enskip
\cd[0.7]{10}{1/6/B,2/4/B,3/10/R,5/7/R,8/9/R}\enskip
\cd[0.7]{10}{1/6/B,2/4/B,3/8/R,5/9/R,7/10/R}\enskip
\cd[0.7]{10}{1/6/B,2/4/B,3/8/R,5/10/R,7/9/R}\enskip
\cd[0.7]{10}{1/6/B,2/4/B,3/9/R,5/8/R,7/10/R}\enskip
\cd[0.7]{10}{1/6/B,2/4/B,3/10/R,5/8/R,7/9/R}\\[4pt]
\cd[0.7]{10}{1/6/B,2/4/B,3/9/R,5/10/R,7/8/R}\enskip
\cd[0.7]{10}{1/6/B,2/4/B,3/10/R,5/9/R,7/8/R}\enskip
\cd[0.7]{10}{1/7/B,2/4/B,3/6/R,5/8/R,9/10/R}\enskip
\cd[0.7]{10}{1/7/B,2/4/B,3/6/R,5/9/R,8/10/R}\enskip
\cd[0.7]{10}{1/7/B,2/4/B,3/6/R,5/10/R,8/9/R}\enskip
\cd[0.7]{10}{1/9/B,2/4/B,3/6/R,5/8/R,7/10/R}\enskip
\cd[0.7]{10}{1/9/B,2/4/B,3/6/R,5/10/R,7/8/R}\enskip
\cd[0.7]{10}{1/8/B,2/4/B,3/9/R,5/6/R,7/10/R}\\[4pt]
\cd[0.7]{10}{1/8/B,2/4/B,3/10/R,5/6/R,7/9/R}\enskip
\cd[0.7]{10}{1/9/B,2/4/B,3/7/R,5/8/R,6/10/R}\enskip
\cd[0.7]{10}{1/9/B,2/4/B,3/7/R,5/10/R,6/8/R}\enskip
\cd[0.7]{10}{1/8/B,2/4/B,3/9/R,5/7/R,6/10/R}\enskip
\cd[0.7]{10}{1/8/B,2/4/B,3/10/R,5/7/R,6/9/R}\enskip
\cd[0.7]{10}{1/8/B,2/4/B,3/9/R,5/10/R,6/7/R}\enskip
\cd[0.7]{10}{1/8/B,2/4/B,3/10/R,5/9/R,6/7/R}\enskip
\cd[0.7]{10}{1/9/B,2/4/B,3/8/R,5/10/R,6/7/R}\\[4pt]
\cd[0.7]{10}{1/5/B,2/6/B,3/4/R,7/8/R,9/10/R}\enskip
\cd[0.7]{10}{1/5/B,2/6/B,3/4/R,7/10/R,8/9/R}\enskip
\cd[0.7]{10}{1/5/B,2/7/B,3/4/R,6/8/R,9/10/R}\enskip
\cd[0.7]{10}{1/5/B,2/7/B,3/4/R,6/10/R,8/9/R}\enskip
\cd[0.7]{10}{1/5/B,2/8/B,3/4/R,6/7/R,9/10/R}\enskip
\cd[0.7]{10}{1/5/B,2/8/B,3/4/R,6/10/R,7/9/R}\enskip
\cd[0.7]{10}{1/6/B,2/7/B,3/4/R,5/8/R,9/10/R}\enskip
\cd[0.7]{10}{1/6/B,2/7/B,3/4/R,5/10/R,8/9/R}\\[4pt]
\cd[0.7]{10}{1/6/B,2/8/B,3/4/R,5/9/R,7/10/R}\enskip
\cd[0.7]{10}{1/8/B,2/6/B,3/4/R,5/9/R,7/10/R}\enskip
\cd[0.7]{10}{1/8/B,2/6/B,3/4/R,5/10/R,7/9/R}\enskip
\cd[0.7]{10}{1/7/B,2/8/B,3/4/R,5/10/R,6/9/R}\enskip
\cd[0.7]{10}{1/5/B,2/7/B,3/6/R,4/9/B,8/10/R}\enskip
\cd[0.7]{10}{1/5/B,2/7/B,3/8/R,4/9/B,6/10/R}\enskip
\cd[0.7]{10}{1/5/B,2/6/B,3/7/R,4/9/R,8/10/R}\enskip
\cd[0.7]{10}{1/5/B,2/6/B,3/10/R,4/7/R,8/9/R}\\[4pt]
\cd[0.7]{10}{1/5/B,2/6/B,3/9/R,4/8/R,7/10/R}\enskip
\cd[0.7]{10}{1/5/B,2/7/B,3/6/R,4/10/R,8/9/R}\enskip
\cd[0.7]{10}{1/5/B,2/8/B,3/6/R,4/10/R,7/9/R}\enskip
\cd[0.7]{10}{1/5/B,2/7/B,3/9/R,4/6/R,8/10/R}\enskip
\cd[0.7]{10}{1/5/B,2/7/B,3/8/R,4/10/R,6/9/R}\enskip
\cd[0.7]{10}{1/5/B,2/7/B,3/9/R,4/8/R,6/10/R}\enskip
\cd[0.7]{10}{1/5/B,2/7/B,3/10/R,4/9/R,6/8/R}\enskip
\cd[0.7]{10}{1/5/B,2/8/B,3/7/R,4/9/R,6/10/R}\\[4pt]
\cd[0.7]{10}{1/5/B,2/8/B,3/10/R,4/7/R,6/9/R}\enskip
\cd[0.7]{10}{1/5/B,2/8/B,3/10/R,4/9/R,6/7/R}\enskip
\cd[0.7]{10}{1/7/B,2/5/B,3/6/R,4/9/R,8/10/R}\enskip
\cd[0.7]{10}{1/7/B,2/5/B,3/6/R,4/10/R,8/9/R}\enskip
\cd[0.7]{10}{1/8/B,2/5/B,3/6/R,4/9/R,7/10/R}\enskip
\cd[0.7]{10}{1/8/B,2/5/B,3/6/R,4/10/R,7/9/R}\enskip
\cd[0.7]{10}{1/7/B,2/5/B,3/9/R,4/6/R,8/10/R}\enskip
\cd[0.7]{10}{1/7/B,2/5/B,3/10/R,4/6/R,8/9/R}\\[4pt]
\cd[0.7]{10}{1/8/B,2/5/B,3/10/R,4/6/R,7/9/R}\enskip
\cd[0.7]{10}{1/8/B,2/5/B,3/7/R,4/9/R,6/10/R}\enskip
\cd[0.7]{10}{1/8/B,2/5/B,3/7/R,4/10/R,6/9/R}\enskip
\cd[0.7]{10}{1/8/B,2/5/B,3/10/R,4/7/R,6/9/R}\enskip
\cd[0.7]{10}{1/6/B,2/7/B,3/5/R,4/9/R,8/10/R}\enskip
\cd[0.7]{10}{1/6/B,2/7/B,3/8/R,4/10/R,5/9/R}\enskip
\cd[0.7]{10}{1/6/B,2/7/B,3/10/R,4/9/R,5/8/R}\enskip
\cd[0.7]{10}{1/3/B,2/4/R,5/7/B,6/8/R,9/10/R}\\[4pt]
\cd[0.7]{10}{1/3/B,2/4/R,5/7/B,6/9/R,8/10/R}\enskip
\cd[0.7]{10}{1/3/B,2/4/R,5/7/B,6/10/R,8/9/R}\enskip
\cd[0.7]{10}{1/3/B,2/4/R,5/9/B,6/7/R,8/10/R}\enskip
\cd[0.7]{10}{1/3/B,2/4/R,5/9/B,6/8/R,7/10/R}\enskip
\cd[0.7]{10}{1/3/B,2/4/R,5/9/B,6/10/R,7/8/R}\enskip
\cd[0.7]{10}{1/3/B,2/4/R,5/6/R,7/8/R,9/10/R}\enskip
\cd[0.7]{10}{1/3/B,2/4/R,5/6/R,7/10/R,8/9/R}\enskip
\cd[0.7]{10}{1/3/B,2/4/R,5/7/R,6/8/B,9/10/R}
\end{center}
\caption{One-cycle chord diagrams for non-orientable genus $g=-5$ ($k=5$), diagrams 241--320 of 473.}
\label{fig:atlas_gm5_4}
\end{figure}

\begin{figure}[ht!]
\begin{center}
\cd[0.7]{10}{1/3/B,2/4/R,5/9/R,6/8/B,7/10/R}\enskip
\cd[0.7]{10}{1/3/B,2/4/R,5/10/R,6/8/B,7/9/R}\enskip
\cd[0.7]{10}{1/3/B,2/4/R,5/7/R,6/9/R,8/10/R}\enskip
\cd[0.7]{10}{1/3/B,2/4/R,5/8/R,6/7/R,9/10/R}\enskip
\cd[0.7]{10}{1/3/B,2/4/R,5/10/R,6/7/R,8/9/R}\enskip
\cd[0.7]{10}{1/3/B,2/4/R,5/8/R,6/10/R,7/9/B}\enskip
\cd[0.7]{10}{1/3/B,2/4/R,5/10/R,6/8/R,7/9/B}\enskip
\cd[0.7]{10}{1/3/B,2/4/R,5/8/R,6/10/R,7/9/R}\\[4pt]
\cd[0.7]{10}{1/3/B,2/4/R,5/9/R,6/8/R,7/10/R}\enskip
\cd[0.7]{10}{1/3/B,2/4/R,5/10/R,6/9/R,7/8/R}\enskip
\cd[0.7]{10}{1/3/B,2/5/R,4/6/R,7/8/R,9/10/R}\enskip
\cd[0.7]{10}{1/3/B,2/5/R,4/6/R,7/10/R,8/9/R}\enskip
\cd[0.7]{10}{1/3/B,2/5/R,4/7/R,6/8/B,9/10/R}\enskip
\cd[0.7]{10}{1/3/B,2/5/R,4/9/R,6/8/B,7/10/R}\enskip
\cd[0.7]{10}{1/3/B,2/5/R,4/7/R,6/9/R,8/10/R}\enskip
\cd[0.7]{10}{1/3/B,2/5/R,4/8/R,6/7/R,9/10/R}\\[4pt]
\cd[0.7]{10}{1/3/B,2/5/R,4/10/R,6/7/R,8/9/R}\enskip
\cd[0.7]{10}{1/3/B,2/5/R,4/8/R,6/10/R,7/9/B}\enskip
\cd[0.7]{10}{1/3/B,2/5/R,4/8/R,6/10/R,7/9/R}\enskip
\cd[0.7]{10}{1/3/B,2/5/R,4/9/R,6/8/R,7/10/R}\enskip
\cd[0.7]{10}{1/3/B,2/5/R,4/10/R,6/9/R,7/8/R}\enskip
\cd[0.7]{10}{1/3/B,2/6/R,4/5/R,7/8/R,9/10/R}\enskip
\cd[0.7]{10}{1/3/B,2/6/R,4/5/R,7/10/R,8/9/R}\enskip
\cd[0.7]{10}{1/3/B,2/7/R,4/5/R,6/9/B,8/10/R}\\[4pt]
\cd[0.7]{10}{1/3/B,2/8/R,4/5/R,6/9/B,7/10/R}\enskip
\cd[0.7]{10}{1/3/B,2/7/R,4/5/R,6/8/R,9/10/R}\enskip
\cd[0.7]{10}{1/3/B,2/7/R,4/5/R,6/10/R,8/9/R}\enskip
\cd[0.7]{10}{1/3/B,2/8/R,4/5/R,6/10/R,7/9/R}\enskip
\cd[0.7]{10}{1/3/B,2/6/R,4/9/R,5/8/B,7/10/R}\enskip
\cd[0.7]{10}{1/3/B,2/6/R,4/10/R,5/8/B,7/9/R}\enskip
\cd[0.7]{10}{1/3/B,2/7/R,4/9/R,5/8/B,6/10/R}\enskip
\cd[0.7]{10}{1/3/B,2/7/R,4/10/R,5/8/B,6/9/R}\\[4pt]
\cd[0.7]{10}{1/3/B,2/8/R,4/9/R,5/7/B,6/10/R}\enskip
\cd[0.7]{10}{1/3/B,2/6/R,4/7/R,5/9/R,8/10/R}\enskip
\cd[0.7]{10}{1/3/B,2/6/R,4/10/R,5/7/R,8/9/R}\enskip
\cd[0.7]{10}{1/3/B,2/6/R,4/8/R,5/10/R,7/9/R}\enskip
\cd[0.7]{10}{1/3/B,2/6/R,4/9/R,5/8/R,7/10/R}\enskip
\cd[0.7]{10}{1/3/B,2/6/R,4/10/R,5/9/R,7/8/R}\enskip
\cd[0.7]{10}{1/3/B,2/7/R,4/6/R,5/9/R,8/10/R}\enskip
\cd[0.7]{10}{1/3/B,2/7/R,4/10/R,5/6/R,8/9/R}\\[4pt]
\cd[0.7]{10}{1/3/B,2/7/R,4/8/R,5/10/R,6/9/R}\enskip
\cd[0.7]{10}{1/3/B,2/7/R,4/10/R,5/9/R,6/8/R}\enskip
\cd[0.7]{10}{1/6/B,2/3/R,4/8/B,5/7/R,9/10/R}\enskip
\cd[0.7]{10}{1/6/B,2/3/R,4/8/B,5/9/R,7/10/R}\enskip
\cd[0.7]{10}{1/6/B,2/3/R,4/9/B,5/8/R,7/10/R}\enskip
\cd[0.7]{10}{1/6/B,2/3/R,4/9/B,5/10/R,7/8/R}\enskip
\cd[0.7]{10}{1/5/B,2/3/R,4/6/R,7/8/R,9/10/R}\enskip
\cd[0.7]{10}{1/5/B,2/3/R,4/6/R,7/10/R,8/9/R}\\[4pt]
\cd[0.7]{10}{1/5/B,2/3/R,4/7/R,6/8/R,9/10/R}\enskip
\cd[0.7]{10}{1/5/B,2/3/R,4/7/R,6/10/R,8/9/R}\enskip
\cd[0.7]{10}{1/5/B,2/3/R,4/8/R,6/7/R,9/10/R}\enskip
\cd[0.7]{10}{1/5/B,2/3/R,4/9/R,6/7/R,8/10/R}\enskip
\cd[0.7]{10}{1/5/B,2/3/R,4/10/R,6/7/R,8/9/R}\enskip
\cd[0.7]{10}{1/5/B,2/3/R,4/8/R,6/10/R,7/9/R}\enskip
\cd[0.7]{10}{1/5/B,2/3/R,4/9/R,6/10/R,7/8/R}\enskip
\cd[0.7]{10}{1/5/B,2/3/R,4/10/R,6/9/R,7/8/R}\\[4pt]
\cd[0.7]{10}{1/7/B,2/3/R,4/5/R,6/9/R,8/10/R}\enskip
\cd[0.7]{10}{1/7/B,2/3/R,4/6/R,5/9/B,8/10/R}\enskip
\cd[0.7]{10}{1/7/B,2/3/R,4/10/R,5/9/B,6/8/R}\enskip
\cd[0.7]{10}{1/6/B,2/3/R,4/7/R,5/9/R,8/10/R}\enskip
\cd[0.7]{10}{1/6/B,2/3/R,4/9/R,5/7/R,8/10/R}\enskip
\cd[0.7]{10}{1/6/B,2/3/R,4/8/R,5/10/R,7/9/R}\enskip
\cd[0.7]{10}{1/6/B,2/3/R,4/10/R,5/8/R,7/9/R}\enskip
\cd[0.7]{10}{1/7/B,2/3/R,4/6/R,5/9/R,8/10/R}\\[4pt]
\cd[0.7]{10}{1/8/B,2/3/R,4/6/R,5/9/R,7/10/R}\enskip
\cd[0.7]{10}{1/8/B,2/3/R,4/6/R,5/10/R,7/9/R}\enskip
\cd[0.7]{10}{1/7/B,2/3/R,4/9/R,5/6/R,8/10/R}\enskip
\cd[0.7]{10}{1/7/B,2/3/R,4/8/R,5/10/R,6/9/R}\enskip
\cd[0.7]{10}{1/7/B,2/3/R,4/9/R,5/8/R,6/10/R}\enskip
\cd[0.7]{10}{1/7/B,2/3/R,4/10/R,5/9/R,6/8/R}\enskip
\cd[0.7]{10}{1/8/B,2/3/R,4/9/R,5/7/R,6/10/R}\enskip
\cd[0.7]{10}{1/8/B,2/3/R,4/10/R,5/7/R,6/9/R}\\[4pt]
\cd[0.7]{10}{1/5/B,2/4/R,3/8/B,6/10/R,7/9/R}\enskip
\cd[0.7]{10}{1/6/B,2/4/R,3/9/B,5/7/R,8/10/R}\enskip
\cd[0.7]{10}{1/6/B,2/4/R,3/8/B,5/10/R,7/9/R}\enskip
\cd[0.7]{10}{1/8/B,2/4/R,3/6/B,5/9/R,7/10/R}\enskip
\cd[0.7]{10}{1/8/B,2/4/R,3/6/B,5/10/R,7/9/R}\enskip
\cd[0.7]{10}{1/7/B,2/4/R,3/9/B,5/8/R,6/10/R}\enskip
\cd[0.7]{10}{1/5/B,2/6/R,3/7/B,4/9/R,8/10/R}\enskip
\cd[0.7]{10}{1/5/B,2/6/R,3/8/B,4/9/R,7/10/R}
\end{center}
\caption{One-cycle chord diagrams for non-orientable genus $g=-5$ ($k=5$), diagrams 321--400 of 473.}
\label{fig:atlas_gm5_5}
\end{figure}

\begin{figure}[ht!]
\begin{center}
\cd[0.7]{10}{1/5/B,2/7/R,3/9/B,4/8/R,6/10/R}\enskip
\cd[0.7]{10}{1/5/B,2/9/R,3/8/B,4/7/R,6/10/R}\enskip
\cd[0.7]{10}{1/6/B,2/5/R,3/8/B,4/9/R,7/10/R}\enskip
\cd[0.7]{10}{1/8/B,2/5/R,3/6/B,4/9/R,7/10/R}\enskip
\cd[0.7]{10}{1/6/B,2/7/R,3/8/B,4/10/R,5/9/R}\enskip
\cd[0.7]{10}{1/4/B,2/5/R,3/7/R,6/9/R,8/10/R}\enskip
\cd[0.7]{10}{1/4/B,2/5/R,3/7/R,6/10/R,8/9/R}\enskip
\cd[0.7]{10}{1/4/B,2/5/R,3/9/R,6/7/R,8/10/R}\\[4pt]
\cd[0.7]{10}{1/4/B,2/5/R,3/9/R,6/8/R,7/10/R}\enskip
\cd[0.7]{10}{1/4/B,2/5/R,3/9/R,6/10/R,7/8/R}\enskip
\cd[0.7]{10}{1/4/B,2/7/R,3/5/R,6/9/R,8/10/R}\enskip
\cd[0.7]{10}{1/4/B,2/7/R,3/5/R,6/10/R,8/9/R}\enskip
\cd[0.7]{10}{1/4/B,2/9/R,3/5/R,6/7/R,8/10/R}\enskip
\cd[0.7]{10}{1/4/B,2/9/R,3/5/R,6/8/R,7/10/R}\enskip
\cd[0.7]{10}{1/4/B,2/9/R,3/5/R,6/10/R,7/8/R}\enskip
\cd[0.7]{10}{1/4/B,2/6/R,3/8/R,5/9/R,7/10/R}\\[4pt]
\cd[0.7]{10}{1/4/B,2/6/R,3/8/R,5/10/R,7/9/R}\enskip
\cd[0.7]{10}{1/4/B,2/8/R,3/6/R,5/9/R,7/10/R}\enskip
\cd[0.7]{10}{1/4/B,2/8/R,3/6/R,5/10/R,7/9/R}\enskip
\cd[0.7]{10}{1/5/B,2/4/R,3/6/R,7/10/R,8/9/R}\enskip
\cd[0.7]{10}{1/5/B,2/4/R,3/7/R,6/10/R,8/9/R}\enskip
\cd[0.7]{10}{1/5/B,2/4/R,3/8/R,6/10/R,7/9/R}\enskip
\cd[0.7]{10}{1/6/B,2/4/R,3/7/R,5/10/R,8/9/R}\enskip
\cd[0.7]{10}{1/6/B,2/4/R,3/10/R,5/7/R,8/9/R}\\[4pt]
\cd[0.7]{10}{1/6/B,2/4/R,3/8/R,5/9/R,7/10/R}\enskip
\cd[0.7]{10}{1/6/B,2/4/R,3/9/R,5/8/R,7/10/R}\enskip
\cd[0.7]{10}{1/7/B,2/4/R,3/8/R,5/10/R,6/9/R}\enskip
\cd[0.7]{10}{1/7/B,2/4/R,3/9/R,5/8/R,6/10/R}\enskip
\cd[0.7]{10}{1/7/B,2/4/R,3/10/R,5/9/R,6/8/R}\enskip
\cd[0.7]{10}{1/7/B,2/8/R,3/4/R,5/10/R,6/9/R}\enskip
\cd[0.7]{10}{1/7/B,2/9/R,3/4/R,5/8/R,6/10/R}\enskip
\cd[0.7]{10}{1/7/B,2/10/R,3/4/R,5/9/R,6/8/R}\\[4pt]
\cd[0.7]{10}{1/5/B,2/6/R,3/9/R,4/8/R,7/10/R}\enskip
\cd[0.7]{10}{1/5/B,2/6/R,3/10/R,4/8/R,7/9/R}\enskip
\cd[0.7]{10}{1/5/B,2/8/R,3/6/R,4/9/R,7/10/R}\enskip
\cd[0.7]{10}{1/5/B,2/9/R,3/8/R,4/6/R,7/10/R}\enskip
\cd[0.7]{10}{1/5/B,2/10/R,3/8/R,4/6/R,7/9/R}\enskip
\cd[0.7]{10}{1/5/B,2/7/R,3/9/R,4/8/R,6/10/R}\enskip
\cd[0.7]{10}{1/5/B,2/9/R,3/8/R,4/7/R,6/10/R}\enskip
\cd[0.7]{10}{1/2/R,3/4/R,5/6/R,7/8/R,9/10/R}\\[4pt]
\cd[0.7]{10}{1/2/R,3/4/R,5/6/R,7/10/R,8/9/R}\enskip
\cd[0.7]{10}{1/2/R,3/4/R,5/7/R,6/9/R,8/10/R}\enskip
\cd[0.7]{10}{1/2/R,3/4/R,5/10/R,6/7/R,8/9/R}\enskip
\cd[0.7]{10}{1/2/R,3/4/R,5/8/R,6/10/R,7/9/R}\enskip
\cd[0.7]{10}{1/2/R,3/4/R,5/10/R,6/9/R,7/8/R}\enskip
\cd[0.7]{10}{1/2/R,3/5/R,4/7/R,6/9/R,8/10/R}\enskip
\cd[0.7]{10}{1/2/R,3/5/R,4/7/R,6/10/R,8/9/R}\enskip
\cd[0.7]{10}{1/2/R,3/5/R,4/9/R,6/7/R,8/10/R}\\[4pt]
\cd[0.7]{10}{1/2/R,3/5/R,4/9/R,6/8/R,7/10/R}\enskip
\cd[0.7]{10}{1/2/R,3/5/R,4/9/R,6/10/R,7/8/R}\enskip
\cd[0.7]{10}{1/2/R,3/6/R,4/5/R,7/10/R,8/9/R}\enskip
\cd[0.7]{10}{1/2/R,3/8/R,4/5/R,6/10/R,7/9/R}\enskip
\cd[0.7]{10}{1/2/R,3/6/R,4/9/R,5/8/R,7/10/R}\enskip
\cd[0.7]{10}{1/2/R,3/6/R,4/10/R,5/8/R,7/9/R}\enskip
\cd[0.7]{10}{1/2/R,3/6/R,4/10/R,5/9/R,7/8/R}\enskip
\cd[0.7]{10}{1/2/R,3/8/R,4/6/R,5/10/R,7/9/R}\\[4pt]
\cd[0.7]{10}{1/2/R,3/10/R,4/6/R,5/8/R,7/9/R}\enskip
\cd[0.7]{10}{1/2/R,3/8/R,4/7/R,5/10/R,6/9/R}\enskip
\cd[0.7]{10}{1/2/R,3/9/R,4/7/R,5/10/R,6/8/R}\enskip
\cd[0.7]{10}{1/2/R,3/10/R,4/7/R,5/9/R,6/8/R}\enskip
\cd[0.7]{10}{1/2/R,3/8/R,4/10/R,5/9/R,6/7/R}\enskip
\cd[0.7]{10}{1/2/R,3/10/R,4/9/R,5/8/R,6/7/R}\enskip
\cd[0.7]{10}{1/3/R,2/5/R,4/8/R,6/10/R,7/9/R}\enskip
\cd[0.7]{10}{1/3/R,2/6/R,4/7/R,5/9/R,8/10/R}\\[4pt]
\cd[0.7]{10}{1/3/R,2/6/R,4/9/R,5/7/R,8/10/R}\enskip
\cd[0.7]{10}{1/3/R,2/6/R,4/8/R,5/10/R,7/9/R}\enskip
\cd[0.7]{10}{1/3/R,2/7/R,4/8/R,5/10/R,6/9/R}\enskip
\cd[0.7]{10}{1/3/R,2/7/R,4/10/R,5/9/R,6/8/R}\enskip
\cd[0.7]{10}{1/4/R,2/6/R,3/8/R,5/9/R,7/10/R}\enskip
\cd[0.7]{10}{1/4/R,2/6/R,3/9/R,5/8/R,7/10/R}\enskip
\cd[0.7]{10}{1/4/R,2/9/R,3/6/R,5/8/R,7/10/R}\enskip
\cd[0.7]{10}{1/4/R,2/8/R,3/7/R,5/10/R,6/9/R}\\[4pt]
\cd[0.7]{10}{1/5/R,2/7/R,3/9/R,4/8/R,6/10/R}
\end{center}
\caption{One-cycle chord diagrams for non-orientable genus $g=-5$ ($k=5$), diagrams 401--473 of 473.}
\label{fig:atlas_gm5_6}
\end{figure}

{\small
\paragraph{Duality, $g=-5$.} 473 diagrams in total: 59 self-dual, 207 dual pairs.
\emph{Self-dual:} 1, 3, 6, 9, 16, 17, 18, 20, 22, 24, 34, 43, 68, 76, 86, 101, 102, 103, 117, 126, 140, 142, 165, 178, 184, 244, 259, 261, 263, 273, 274, 277, 278, 284, 288, 297, 298, 306, 312, 314, 318, 325, 330, 342, 343, 355, 371, 374, 385, 423, 440, 441, 443, 445, 451, 453, 458, 462, 471.
\emph{Dual pairs:} (2,19), (4,23), (5,13), (7,15), (8,12), (10,14), (11,21), (25,48), (26,253), (27,30), (28,354), (29,84), (31,257), (32,55), (33,247), (35,345), (36,356), (37,66), (38,351), (39,79), (40,270), (41,46), (42,255), (44,119), (45,349), (47,96), (49,262), (50,266), (51,339), (52,58), (53,192), (54,93), (56,361), (57,250), (59,258), (60,113), (61,99), (62,334), (63,176), (64,265), (65,73), (67,90), (69,326), (70,72), (71,338), (74,264), (75,271), (77,359), (78,83), (80,353), (81,197), (82,110), (85,344), (87,321), (88,337), (89,341), (91,114), (92,313), (94,347), (95,331), (97,243), (98,185), (100,332), (104,249), (105,252), (106,360), (107,336), (108,120), (109,183), (111,180), (112,172), (115,267), (116,217), (118,209), (121,246), (122,346), (123,168), (124,362), (125,438), (127,467), (128,166), (129,389), (130,146), (131,388), (132,418), (133,469), (134,218), (135,416), (136,431), (137,204), (138,398), (139,282), (141,406), (143,439), (144,289), (145,152), (147,468), (148,190), (149,234), (150,473), (151,390), (153,162), (154,295), (155,432), (156,304), (157,231), (158,219), (159,430), (160,206), (161,200), (163,414), (164,235), (167,329), (169,179), (170,300), (171,316), (173,201), (174,223), (175,335), (177,269), (181,203), (182,268), (186,357), (187,328), (188,444), (189,460), (191,461), (193,310), (194,240), (195,377), (196,420), (198,205), (199,375), (202,409), (207,411), (208,340), (210,256), (211,323), (212,283), (213,392), (214,424), (215,293), (216,248), (220,412), (221,454), (222,449), (224,459), (225,447), (226,450), (227,369), (228,383), (229,442), (230,463), (232,287), (233,429), (236,413), (237,446), (238,457), (239,241), (242,320), (245,260), (251,352), (254,350), (272,333), (275,370), (276,372), (279,455), (280,452), (281,415), (285,291), (286,394), (290,381), (292,400), (294,427), (296,305), (299,303), (301,464), (302,425), (307,435), (308,470), (309,405), (311,472), (315,317), (319,324), (322,327), (348,358), (363,387), (364,386), (365,456), (366,448), (367,373), (368,376), (378,379), (380,410), (382,408), (384,421), (391,407), (393,422), (395,399), (396,417), (397,401), (402,437), (403,465), (404,466), (419,428), (426,434), (433,436).
}


\clearpage 



\section*{Tool use declaration}

The enumeration and duality computations were performed with Python programs
\texttt{one\_cycle\_chord\_diagrams.py} and \texttt{general\_codes.py}
\cite{zenodo22298943, zenodo23084341}; the figures were typeset with \LaTeX\ and TikZ.

Generative AI (Claude, Anthropic) was used as follows: (i) writing and debugging
parts of the Python code, including the scripts that generate the catalog of
chord diagrams in the Appendix (conversion from the computational code to the
drawing code, canonical forms, dual and self-dual diagrams); (ii) generating the
\LaTeX/TikZ source of the Appendix figures from the program output;
(iii) drafting and editing the English text of the description of the algorithms
and of some passages of the Introduction and Conclusion, and suggesting
arXiv subject classes.

All mathematical statements, definitions and results come from the authors and
from the cited literature. The authors checked the program output independently
(for example, by verifying that taking the dual twice returns the original
diagram and that the counts of diagrams are consistent with the known ones),
reviewed and edited all AI-generated text and code, and take full responsibility
for the content of the paper.

\end{document}